\input ./style/arxiv-general.cfg
\documentclass[aap,MSNbibl,dvips]{arximspdf}
\makeatletter
   \@ifpackageloaded{graphicx}{}{\usepackage{graphicx}}
\makeatother
\doi{10.1214/14-AAP1083}% Updated by VTEXPTS2LaTeX.exe, 06.01.2015
\volume{25}
\issue{6}
\pubyear{2015}
\firstpage{3592}
\lastpage{3623}
\docsubty{FLA}

\makeatletter
\renewcommand{\P}{\mathbf{P}}
\newproclaim{exa}{Example}
\newtheorem{prop}[thm]{Proposition}
\newtheorem{lem}[thm]{Lemma}
\newtheorem{cor}[thm]{Corollary}
\newproclaim{rem}{Remark}
\makeatother

\begin{document}
\begin{frontmatter}

\title{Stability of adversarial Markov chains, with an application to
adaptive MCMC algorithms}
\runtitle{Adversarial Markov chains}

\begin{aug}
% Corresponding author: Jeffrey Rosenthal - jeff@math.toronto.edu% Updated by VTEXPTS2LaTeX.exe, 09.01.2015 12:39
%Updated by VTEXPTS2LaTeX.exe, 06.01.2015 13:40
\author[A]{\fnms{Radu V.}~\snm{Craiu}\thanksref{m1}},
\author[B]{\fnms{Lawrence}~\snm{Gray}\thanksref{m2}},
\author[C]{\fnms{Krzysztof}~\snm{{\L}atuszy\'{n}ski}\thanksref{m3}},
\author[D]{\fnms{Neal}~\snm{Madras}\thanksref{m4}},
\author[C]{\fnms{Gareth O.}~\snm{Roberts}\thanksref{m3}}
\and
\author[A]{\fnms{Jeffrey S.}~\snm{Rosenthal}\corref{}\thanksref{m1}\ead[label=e1]{jeff@math.toronto.edu}}
\runauthor{R. V. Craiu et al.}
\affiliation{University of Toronto\thanksmark{m1}, University of Minnesota\thanksmark{m2},
University of Warwick\thanksmark{m3} and  York~University\thanksmark{m4}}
%\dedicated{}
\address[A]{R. V. Craiu \\
J. S. Rosenthal\\
Department of Statistics\\
University of Toronto\\
100 Street George Street\\
Toronto, Ontario M5S 3G3\\
Canada\\
\printead{e1}}
\address[B]{L. Gray\\
Department of Mathematics\\
University of Minnesota\\
206 Church Street SE\\
Minneapolis, Minnesota 55455-0488\\
USA}
\address[C]{K. {\L}atuszy\'{n}ski\\
G. O. Roberts\\
Department of Statistics\\
University of Warwick\\
Coventry CV4 7AL\\
United Kingdom}
\address[D]{N. Madras\\
Department of Mathematics\hspace*{28pt}\\
\quad and Statistics\\
York University\\
4700 Keele St.\\
Toronto, Ontario M3J 1P3\\
Canada}
\end{aug}

% HISTORY:
%
\received{\smonth{3} \syear{2014}}% Updated by VTEXPTS2LaTeX.exe,
%06.01.2015 13:40
%
\revised{\smonth{8} \syear{2014}}% Updated by VTEXPTS2LaTeX.exe,
%06.01.2015 13:40

% ABSTRACT
%
\begin{abstract}
We consider whether ergodic Markov chains with bounded step size
remain bounded in probability when their transitions are modified by
an adversary on a bounded subset. We provide counterexamples to show
that the answer is no in general, and prove theorems to show that the
answer is yes under various additional assumptions. We then use our
results to prove convergence of various adaptive Markov chain Monte
Carlo algorithms.
\end{abstract}

% KEYWORDS
\end{frontmatter}

%s1 #&#
\section{Introduction}

This paper considers whether bounded modifications of stable Markov
chains remain stable. Specifically, we let $P$ be a
fixed time-homogeneous ergodic Markov chain kernel with bounded step
size, and let $\{X_n\}$ be a stochastic process which follows the
transition probabilities $P$ except on a bounded subset $K$ where an
``adversary'' can make arbitrary bounded jumps. Under what conditions
must such a process $\{X_n\}$ be bounded in probability?

One might think that such
boundedness would follow easily, at least under mild regularity and
continuity assumptions, that is, that modifying a stable continuous Markov
chain inside a bounded set $K$ couldn't possibly lead to unstable
behavior out in the tails. In fact the situation is rather more
subtle, as we explore herein. We will provide counterexamples to
show that boundedness may fail even for well-behaved continuous
chains. We will then show that under various additional conditions,
including bounds on transition probabilities and/or small set
assumptions and/or geometric ergodicity, such boundedness does hold.

The specific question considered here appears to be new, though it is
somewhat reminiscent of previous bounds on non-Markovian stochastic
processes such as those related to \textit{adversarial queuing
theory} \cite{hajek,pemantle,borodin}.
We present our formal setup in Section~\ref{sec-setup}, our
main results in Section~\ref{sec-results} and some counterexamples
in Section~\ref{sec-counterex}. Our results are then
proven in Sections~\ref{sec-cemproof} through \ref{sec-geomproof}.

In Section~\ref{sec-mcmc}, we turn our attention to
adaptive Markov chain Monte Carlo (MCMC) algorithms.
MCMC proceeds by running
a Markov chain long enough to approximately converge to
its stationary distribution and thus provide useful samples.
Adaptive MCMC algorithms attempt to improve on MCMC by
modifying the Markov chain transitions as they run, but this destroys
the Markov property and makes convergence to stationarity notoriously
difficult to prove.
We use our main results herein to establish general conditions for
convergence of certain adaptive MCMC algorithms
(Theorem~\ref{mcmcthm}). We then apply this result to a simple
but useful adaptive MCMC algorithm (Proposition~\ref{BAMprop}),
and also to a detailed statistical application involving a probit
model for lupus patient data (Section~\ref{sec-statexample}).
For details and references about adaptive MCMC algorithms,
see Section~\ref{sec-mcmc}.

%s2 #&#
\section{Formal setup and assumptions}
\label{sec-setup}

% Let $(\X,\F)$ be a metricised measure space, and

Let $\mathcal{X}$ be a nonempty general (i.e., possibly uncountable) state
space, on which is defined a metric $\eta$, which gives rise to a
corresponding Borel $\sigma$-algebra $\mathcal{F}$.
Assume that $\mathcal{X}$ contains some specified ``origin'' point
$\mathbf{0}\in\mathcal{X}$.
(In our examples and applications, $\mathcal{X}$ will usually be a
subset of
$\mathbf{R}^d$ with the usual Euclidean metric.)
Let $P$ be the transition probability kernel for a fixed
time-homogeneous Markov chain on $\mathcal{X}$.
Assume that $P$ is Harris ergodic with
stationary probability distribution $\pi$, so that
%
%e1 #&#
\begin{equation}
\label{eqn-ergodic} \lim_{n\to\infty} \bigl\|P^n(x,\cdot) - \pi\bigr\| :=
\lim_{n\to\infty} \sup_{A \in\mathcal{F}} \bigl|P^n(x,A) -
\pi(A) \bigr| = 0,\qquad x \in\mathcal{X}.
\end{equation}
We assume, to relate the Markov chain to the geometry of $\mathcal{X}$,
that there is a constant $D<\infty$ such that
$P$ never moves more than a distance $D$, that is, such that
%
%e2 #&#
\begin{equation}
\label{eqn-boundedjumps} P \bigl( x, \bigl\{y \in\mathcal{X}\dvtx\eta(x,y) \le D \bigr\}
\bigr) = 1,\qquad x \in\mathcal{X}.
\end{equation}
Let $K\in\mathcal{F}$ be a fixed \textit{bounded} nonempty subset of
$\mathcal{X}$,
and for $r>0$ let $K_r$ be the set of all states within a
distance $r$ of $K$ (so each $K_r$ is also bounded).

In terms of these ingredients, we define our
``adversarial Markov chain'' process $\{X_n\}$ as follows.
It begins with $X_0 = x_0$ for some specific initial state $x_0$;
for simplicity (see the proof of Lemma~\ref{lemma-cembound})
we assume that $x_0 \in K$.
Whenever the process is outside of $K$,
it moves according to the Markov transition probabilities $P$,
that~is,
%
%e3 #&#
\begin{equation}
\label{XfollowsP}
\hspace*{6pt}\P(X_{n+1} \in A | X_0,X_1,
\ldots,X_n) = P(X_n,A),\qquad n \ge0, A\in\mathcal{F},
X_n \notin K.
\end{equation}
When the process is inside of $K$, it can move \textit{arbitrarily},
according to an adversary's wishes, perhaps depending
on the time $n$ and/or the chain's history in a nonanticipatory
manner (i.e., adapted to $\{X_n\}$; see also Example~\ref{exa3} below),
subject only to measurability [i.e.,
$\P(X_{n+1} \in A | X_0,X_1,\ldots,X_n)$ must be well defined
for all $n \ge0$ and $A\in\mathcal{F}$], and to
the restriction that it can't move more than a distance $D$
at each iteration---or more specifically
that from $K$, it can only move to points within $K_D$.
In summary, $\{X_n\}$ is a stochastic process which is
``mostly'' a Markov chain following the
transition probabilities~$P$, except that it is modified
by an adversary when it is within the bounded subset $K$.

We are interested in conditions guaranteeing that this process
$\{X_n\}$ will be bounded in probability, that is, will be tight, that is,
will satisfy that
%
%e4 #&#
\begin{equation}
\label{boundedinprob} \lim_{L\to\infty} \sup_{n\in\mathbf{N}} \P\bigl(
\eta(X_n,\mathbf{0}) > L | X_0=x_0\bigr) =
0.
\end{equation}

%s3 #&#
\section{Main results}
\label{sec-results}

We now consider various conditions under which (\ref{boundedinprob})
will or will not hold. For application of our results to the
verification of adaptive MCMC algorithms, see Section~\ref{sec-mcmc}
below.

%s3.1 #&#
\subsection{First results}

We first note that such boundedness is guaranteed in the
\textit{absence} of an adversary:

%pr1 #&#
\begin{prop}
\label{noadprop}
In the setup of Section~\ref{sec-setup},
suppose $\{X_n\}$ always follows the transitions $P$ (including when
it is within $K$, i.e., there is no adversary).
Then (\ref{boundedinprob}) holds.
\end{prop}

% \noindent Proposition~\ref{noadprop} is proved in
% Section~\ref{sec-prop1proof} below.

Indeed, Proposition~\ref{noadprop} follows immediately since
if $P$ is Harris ergodic as
in~(\ref{eqn-ergodic}), then it converges in distribution, so it must
be tight and hence satisfy~(\ref{boundedinprob}).
[In fact, even if $P$ is just assumed to be
$\phi$-irreducible with period $d \ge1$ and stationary
probability distribution $\pi$, then this argument can be applied
separately to each of the sequences $\{X_{dn+j}\}_{n=0}^\infty$
for $j=0,1,\ldots,d-1$ to again conclude~(\ref{boundedinprob}).]

Boundedness also holds for a lattice like $\mathbf{Z}^d$, or more generally
if the state space $\mathcal{X}$ is topologically discrete (i.e., countable
and such that each state $x$ is topologically isolated and hence
open in $\mathcal{X}$). In this case, bounded subsets like $K_{2D}$
must be
\textit{finite}, and the result holds without any further assumptions:

%pr2 #&#
\begin{prop}
\label{discreteprop}
In the setup of Section~\ref{sec-setup},
suppose $P$ is an irreducible positive-recurrent Markov chain with
stationary probability distribution $\pi$ on a countable state
space $\mathcal{X}$ such that $K_{2D}$ is finite.
Then (\ref{boundedinprob}) holds.
\end{prop}

Proposition~\ref{discreteprop} is proved in
Section~\ref{sec-cemproof} below.

However, (\ref{boundedinprob}) does not hold in general, not even
under a strong continuity assumption:

%pr3 #&#
\begin{prop}
\label{counterexprop}
There exist adversarial Markov chain examples following
the setup of Section~\ref{sec-setup}, on state spaces which are
countable subsets of $\mathbf{R}^2$, which fail to
satisfy~(\ref{boundedinprob}), even under
the strong continuity condition that $\mathcal{X}$ is closed and
%
%e5 #&#
\begin{eqnarray}
&&\forall x\in\mathcal{X}, \forall \epsilon>0, \exists \delta>0
\quad\mbox{s.t.}\quad \bigl\|P(y,\cdot)-P(x,\cdot)\bigr\| < \epsilon \nonumber
\\[-8pt]
\label{contcond}
\\[-8pt]
\eqntext{\mbox{whenever }
\eta(x,y)<\delta.}
\end{eqnarray}
\end{prop}

Proposition~\ref{counterexprop} is proved in
Section~\ref{sec-counterex} below, using two different counterexamples.

Proposition~\ref{counterexprop} says that the adversarial process
$\{X_n\}$ may not be bounded in probability, even if we assume a
strong continuity condition on $P$. Hence, additional assumptions
are required, as we consider next.

%re1 #&#
\begin{rem}
The counterexamples in Proposition~\ref{counterexprop} are discrete
Markov chains in the sense that their state spaces are countable.
However, their state spaces~$\mathcal{X}$ are not \textit{topologically}
discrete, since they contain accumulation points, and in particular
sets like $K_{2D}$ are not finite there, so there is no contradiction
with Proposition~\ref{discreteprop}.
\end{rem}

%s3.2 #&#
\subsection{A result using expected hitting times}

We now consider two new assumptions. The first provides an upper
bound on the Markov chain transitions out of~$K_D$:

\begin{longlist}[(A1)]
\item[(A1)] There is $M<\infty$, and a probability measure $\mu_*$
concentrated
on $K_{2D} \setminus K_D$, such that
$P(x,dz) \le M \mu_*(dz)$ for all $x \in K_D \setminus K$ and $z \in
K_{2D} \setminus K_D$.
\end{longlist}

Note that in (A1) we always have $z \neq x$,
which is helpful when considering, for example, Metropolis algorithms which
have positive probability of not moving.
Choices of $\mu_*$ in (A1) might include
$\operatorname{Uniform} (K_{2D} \setminus K_D)$,
or $\pi|_{K_{2D} \setminus K_D}$.
The second assumption bounds an expected hitting time:

\begin{longlist}[(A2)]
\item[(A2)] The expected time for a Markov chain following the
transitions $P$
to reach the subset $K_D$, when started from the distribution $\mu_*$
in (A1), is finite.
\end{longlist}

In terms of these two assumptions, we have:

%th4 #&#
\begin{thm}
\label{thm-cemetery}
In the setup of Section~\ref{sec-setup},
if \textup{(A1)} and \textup{(A2)} hold for the same~$\mu_*$,
then (\ref{boundedinprob}) holds; that is,
$\{X_n\}$ is bounded in probability.
\end{thm}

Theorem~\ref{thm-cemetery} is proved in
Section~\ref{sec-cemproof} below.

%s3.3 #&#
\subsection{A result assuming a small set condition}

Condition (A2), that the hitting time of $K_D$
has finite expectation, may be difficult to verify directly. As an
alternative, we consider a different assumption:

\begin{longlist}[(A3)]
\item[(A3)] The set $K_{2D} \setminus K_D$ is small for $P$; that is, there
is some probability measure~$\nu_*$ on $\mathcal{X}$, and
some $\epsilon>0$, and some $n_0\in\mathbf{N}$, such that
$P^{n_0}(x,A) \ge \epsilon \nu_*(A)$
for all states $x \in K_{2D} \setminus K_D$ and all subsets $A \in
\mathcal{F}$.
\end{longlist}

We then have:
%
%th5 #&#
\begin{thm}
\label{thm-pismall}
In the setup of Section~\ref{sec-setup},
if \textup{(A1)} and \textup{(A3)} hold where
either \textup{(a)}~$\nu_*=\mu_*$, or \textup{(b)} $P$ is reversible and
$\mu_*=\pi|_{K_{2D} \setminus K_D}$,
then (\ref{boundedinprob}) holds; that is,
$\{X_n\}$ is bounded in probability.
\end{thm}

Theorem~\ref{thm-pismall} is proved in
Section~\ref{sec-smallproof} below.

Assumption~(A3) is often straightforward
to verify. For example:

%pr6 #&#
\begin{prop} \label{unifsmallprop}
Suppose $\mathcal{X}$ is an open subset of $\mathbf{R}^d$ which
contains a~bounded
rectangle $J$ which contains $K_{2D} \setminus K_D$.
Suppose there are $\delta>0$ and $\epsilon>0$ such that
%
%e6 #&#
\begin{equation}
\label{deltepeqn} P(x,dy) \ge \epsilon \operatorname{Leb}(dy) \qquad\mbox{whenever }
x,y\in J \mbox{ with } |y-x|<\delta,
\end{equation}
where $\operatorname{Leb}$ is Lebesgue measure on $\mathbf{R}^d$.
Then \textup{(A3)} holds with $\nu_* = \break \operatorname{Uniform}(K_{2D}
\setminus K_D)$.
\end{prop}

Proposition~\ref{unifsmallprop} is proved
in Section~\ref{sec-unifsmallproof} below.

%s3.4 #&#
\subsection{A result assuming geometric ergodicity}

Assumption (A3) can be verified for various Markov chains,
as we will see below. However, its verification will sometimes be
difficult. An alternative approach is to consider \textit{geometric
ergodicity}, as follows (see, e.g., \cite{MT} for context):

\begin{longlist}[(A4)]
\item[(A4)] The Markov chain transition kernel
$P$ is geometrically ergodic; that is, there is $\rho<1$
and a $\pi$-a.e. finite measurable function $\xi\dvtx\mathcal{X}\to
[1,\infty]$ such that
$\|P^n(x,\cdot) - \pi\| \le\xi(x) \rho^n$ for $n\in\mathbf{N}$
and $x\in\mathcal{X}$.
\end{longlist}

We also require a slightly different version of
(A1):

\begin{longlist}[(A5)]
\item[(A5)]
There is $M<\infty$ such that
$P(x,dz) \le M \pi(dz)$ for all
$x\in K_D$ and $z\in K_{2D}$.
\end{longlist}

[Of course, (A5) holds trivially for $z \notin K_{2D}$,
since then $P(x,dz)=0$.]
We then have:

%th7 #&#
\begin{thm}
\label{thm-geomerg}
In the setup of Section~\ref{sec-setup},
if \textup{(A4)} and \textup{(A5)} hold,
then (\ref{boundedinprob}) holds; that is,
$\{X_n\}$ is bounded in probability.
\end{thm}

Theorem~\ref{thm-geomerg} is proved in
Section~\ref{sec-geomproof} below.

% \subsection{A short proof of Proposition~\ref{noadprop}}
% \labelss{sec-prop1proof}
%
% We end this section with the (easy) proof of Proposition~
%\ref{noadprop}.
%
% Let $x_0\in K$ be the initial state, let $B_L = \{x\in\X: \eta(x,
%\bzero) \ge
% L\}$ and let $\epsilon>0$. Then by equation~\eqref{eqn-ergodic},
% we can find $N$ large enough that
% $\|P^n(x_0,\cdot)-\pi\| < \epsilon/2$ for all $n \ge N$.
% Then choose $L$ large enough
% that $\pi(B_L) < \epsilon/2$ and also $L > |x_0| + ND$.
% For $n<N$, it follows
% from equation~\eqref{eqn-boundedjumps} that $P^n(x_0,B_L)=0$.
% For $n \ge N$, we have $P^n(x_0,B_L) \le\pi(B_L) + \|P^n(x_0,\cdot)-
%\pi\|
% < \epsilon/2 + \epsilon/2 = \epsilon$. Hence, $\sup_n P^n(x_0,B_L)
% \le\epsilon$, i.e.\ $\sup_n \P(X_n \in B_L | X_0=x_0) \le\epsilon$.
% Since $\epsilon>0$ was arbitrary, \eqref{boundedinprob} follows.
%
% (Alternatively, if $P$ is Harris ergodic as
% in~\eqref{eqn-ergodic}, then it converges in distribution, so it must
% be tight and hence satisfy~\eqref{boundedinprob}. In fact, even if
% $P$ is just assumed to be
% $\phi$-irreducible with period $d \ge1$ and stationary
% probability distribution $\pi$, then this argument can be applied
% separately to each of the sequences $\{X_{dn+j}\}_{n=0}^\infty$
% for $j=0,1,\ldots,d-1$ to again conclude~\eqref{boundedinprob}.)
% \qed

%s4 #&#
\section{Counterexamples to prove Proposition~\texorpdfstring{\protect\ref{counterexprop}}{3}}
\label{sec-counterex}

We next present two counterexamples to illustrate that with the
setup and assumptions of Section~\ref{sec-setup}, the bounded in
probability property (\ref{boundedinprob}) might fail. Each example
has a state space $\mathcal{X}$ which is a countable subset of
$\mathbf{R}^2$ with the
usual Euclidean metric $\eta(x,y) := |y-x|$. In Example~\ref{exa1},
$\mathcal{X}$
is not closed, and (\ref{contcond}) does not hold; this is remedied
in Example~\ref{exa2}.

%ex1 #&#
\begin{exa}\label{exa1}
Let $\mathcal{X}= \{(\frac{1}{i}, j) \dvtx i\in\mathbf{N},
j=0,1,\ldots\,\}$ be the state space.
That is, $\mathcal{X}= \bigcup_{i\in\mathbf{N}} \mathcal{X}_i$
where each
$\mathcal{X}_i \equiv\{(\frac{1}{i}, j)\}_{j=0,1,\ldots}$ is a
different ``column.''
Let $\pi(\frac{1}{i}, j) = 2^{-i} (\frac{1}{i}) (1-\frac{1}{i})^j$,
so that
$\pi$ restricted to each $\mathcal{X}_i$ is a geometric
distribution with mean $i$.
Let $K = \{(\frac{1}{i}, 0)\}$ consist of the bottom element of each
column; see Figure~\ref{ex1diag}.

%
%f1 #&#
\begin{figure}

\includegraphics{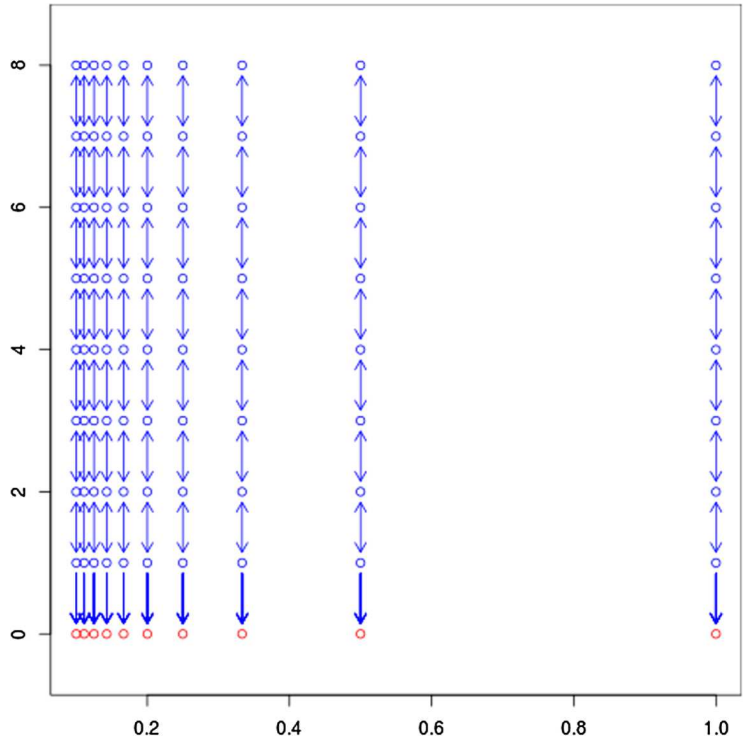}

\caption{Part of the state space in Example~\protect\ref{exa1}.}
\label{ex1diag}
\end{figure}

Let the Markov chain $P$ proceed, outside of $K$, by doing a simple
$\pm1$
Metropolis algorithm up and down its current column $\mathcal{X}_i$ to be
reversible\vspace*{1pt} with respect to $\pi$. That is,\vspace*{1pt} for $j
\ge1$, $P((\frac{1}{i}, j), (\frac{1}{i}, j-1))=\frac{1}{2}$ and
$P((\frac{1}{i}, j),\break  (\frac{1}{i}, j+1))=\frac{1}{2}(1-\frac{1}{i})$,
and the leftovers $P((\frac{1}{i}, j), (\frac{1}{i}, j))
= 1 - P((\frac{1}{i}, j), (\frac{1}{i}, j-1)) - P((\frac{1}{i},
j),(\frac{1}{i}, j+1))$.
Intuitively, the larger the column number $i$, the higher is the
conditional mean of $\pi$ on $\mathcal{X}_i$, so the higher the chain will
tend to move within $\mathcal{X}_i$, and the longer it will take to
return to $K$.

Inside of $K$, choose any appropriate transitions to make the chain
irreducible and reversible with respect to $\pi$; for example, choose
$P ((\frac{1}{i}, 0),(\frac{1}{i}, 1) )=\frac{1}{2}(1-\frac{1}{i})$
and $P ((\frac{1}{i}, 0),(\frac{1}{i-1}, 0) )=1/4$
(for $i>1$ only, otherwise~0),
and $P ((\frac{1}{i}, 0),(\frac{1}{i+1}, 0) ) = i/8(i+1)$,
and the leftovers $P ((\frac{1}{i}, 0),(\frac{1}{i}, 0) )
= 1 - \break  P ((\frac{1}{i}, 0),(\frac{1}{i+1}, 0))
- P ((\frac{1}{i}, 0), (\frac{1}{i-1}, 0))
- P ((\frac{1}{i}, 0), (\frac{1}{i}, 1))$.

Let the adversary proceed within $K$ as follows.
If $X_n \in K$, then $X_{n+1} = (\frac{1}{n}, 1)$.
That is, the chain moves from $K$ to higher and
higher column numbers as time goes on.

With these specifications,
$K$ is bounded, and the process $\{X_n\}$
never moves more than a distance $D=1$, so the setup of
Section~\ref{sec-setup} is satisfied.
However, the process $\{X_n\}$ will,
over time, move closer and closer to~0 in the $x$-direction,
and will then tend to climb higher and higher in the $y$-direction.
More formally, write $X_{n,1}$ and $X_{n,2}$ for the
$x$-coordinate and $y$-coordinate of $X_n$. Then
given any \mbox{$L<\infty$},
choose $m\in\mathbf{N}$ such that the median of a mean-$m$ Geometric random
variable, $\lceil-1/\log_2(1-{1 \over m}) \rceil$, is at least $L$.
Then let $\tau= \inf\{ n \dvtx X_{n,1} \le{1 \over m}\}$. Then
after time $\tau$, the $y$-coordinate of $X_n$ will be stochastically
larger than a usual $\pm1$ Metropolis algorithm for a Geometric
distribution with mean $m$. Hence,
$\liminf_{n\to\infty} \P(X_{n,2} \ge L)$ will be at least as large as
the probability that a mean-$m$ Geometric random variable will
be $\ge L$. This probability is at least $\frac{1}{2}$.
It follows that $\{X_{n,2}\}$, and hence also $\{X_n\}$, are not
bounded in probability, that is, that (\ref{boundedinprob}) does not hold.

(Alternatively, the adversary could proceed within $K$ by moving from
$(\frac{1}{i}, 0)$ to either $(\frac{1}{i}, 1)$ with probability
$\frac{1}{2}(1-\frac{1}{i})$,
or to $(\frac{1}{i+1}, 0)$ with probability $(1+\frac{1}{i})/4$,
or to $(\frac{1}{i-1}, 0)$
with probability $(1+\frac{1}{i})/4$ [for $i>1$ only, otherwise~0],
or remain at $(\frac{1}{i}, 0)$ with the leftover probability.
This would
make the process $\{X_n\}$ be time-homogeneous Markov and reversible with
respect to the infinite measure $\overline\pi$ defined by
$\overline\pi(\frac{1}{i}, j) = (\frac{1}{i}) (1-\frac{1}{i})^j$. Then
$\{X_n\}$ will therefore be null recurrent.
Hence, again, (\ref{boundedinprob}) will not hold.)
\end{exa}

Now, in the above example, the state space $\mathcal{X}$ is not closed.
One could easily ``extend'' the example to include $\{(0,j)\dvtx
j \in\mathbf{N}\}$ and thus make $\mathcal{X}$ closed. However, this cannot
be done in a continuous way; that is, there is no way to satisfy
(\ref{contcond}) in this example. This might lead one to suspect
that a continuity condition such as (\ref{contcond}) suffices to
guarantee (\ref{boundedinprob}). However, that is not the
case, as the following example shows:

%ex2 #&#
\begin{exa}\label{exa2}
Our state space $\mathcal{X}$ will be another countable subset of
$\mathbf{R}^2$,
defined as follows. Let $O=(0,0)$ be the origin. Let $S_0
= \{(i,0) \dvtx i \in\mathbf{N}\}$.
Let $\{\beta_k\}_{k=1}^\infty$ be an increasing sequence of
integers with $\beta_k > k$ to be specified later.
For $k\in\mathbf{N}$, let $S_k$ consist of the $k$
points $(i, {i \over k})$ for $i=1,2,\ldots,k$, together with $\beta_k-1$
additional points equally spaced on the line segment
from $(k,1)$ to the $y$-axis point $(0, \beta_k)$.
Finally, let $\mathcal{Y}=\{(0,i) \dvtx i \in\mathbf{N}\}$ be the
positive integer $y$-axis.
Then $\mathcal{X}= O \cup\mathcal{Y}\cup\bigcup_{k=0}^\infty S_k$;
see Figure~\ref{ex2diag}.

%
%f2 #&#
\begin{figure}

\includegraphics{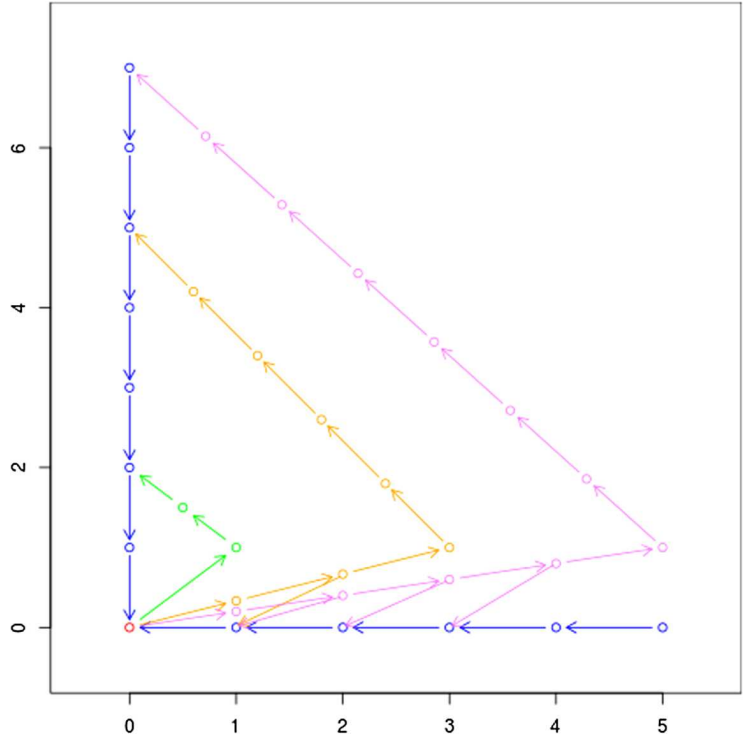}

\caption{Part of the state space in Example~\protect\ref{exa2}, including
$O$ (origin), and $\mathcal{Y}$ ($y$-axis),
and $S_0$ ($x$-axis),
and $S_1$
with $\beta_1=2$ [through $(1,1)$],
and $S_3$ with $\beta_3=5$
[through $(3,1)$],
and $S_5$ with $\beta_5=7$ [through $(5,1)$].}\vspace*{-6pt}
\label{ex2diag}
\end{figure}

Define transitions $P$ on $\mathcal{X}$ as follows.
On $S_0$, we have $P ( (i,0), (i-1, 0) ) = 1$; that is, it always
moves toward the origin.
Similarly, on $\mathcal{Y}$, we have $P ( (0,i), (0,i-1) ) = 1$;
that is, it again always moves toward the origin.\vspace*{1pt}
On the first $k-1$ points of $S_k$, we have
$P ( (i, {i \over k}), (i+1, {i+1 \over k}) ) = {i \over k}$,
and $P ( (i, {i \over k}), (i-1,0) ) = 1 - {i \over k}$; that is, it
either continues upwards on $S_k$, or
moves toward the origin on $S_0$.
On the remaining points of $S_k$, with probability~1
it moves one additional point along $S_k$'s path
toward $(0, \beta_k)$.
The chain's step sizes are thus all bounded above by, for example,
$D=\sqrt{2}$.

Note that these transition probabilities are continuous
in a very strong sense: if $x_n \to x$ (which can only happen for $x
\in S_0$ or for a.a. constant sequences), then $P(x_n,y) \to P(x,y)$ for
all $y\in\mathcal{X}$, and in particular $\|P(x_n,\cdot) - P(x,\cdot
)\| \to0$.
So, (\ref{contcond}) is satisfied.

Note also that if this chain is started at $(1, {1 \over k})$, then
it has probability $\prod_{i=1}^k ({i \over k}) > 0$ of continuing along
$S_k$ all the way to $(k,1)$, in which case it will take a total
of $k+2\beta_k$ iterations to return to $O$. Otherwise,
for $1 \le j \le k-1$, it takes $2j-1$ iterations with probability
$ ( \prod_{i=1}^{j-1} {i \over k} ) ( 1-{j \over k} )$.
Thus, if $r_k = \mathbf{E}(\tau_O | X_0 = (1, {1 \over k}))$ is the
expected return time to $O$ from $(1, {1 \over k})$, then
\[
r_k = (k+2\beta_k) \Biggl( \prod
_{i=1}^k {i \over k} \Biggr) + \sum
_{j=1}^{k-1} ( 2j-1 ) \Biggl( \prod
_{i=1}^{j-1} {i \over k} \Biggr) \biggl( 1-
{j \over k} \biggr).
\]
In particular, by letting $\beta_k$ grow sufficiently quickly, we can
make $r_k$ grow as quickly as desired.

Finally, we specify that from $O$, for $k\in\mathbf{N}$ the Markov chain
moves to $(1, {1 \over k})$ with probability $a_k$, for some positive
numbers $a_k$ summing to~1 to be specified later.

Meanwhile, the adversary's compact set is given by the single state
$K=\{O\}$. From $O$, the adversary proceeds simply by moving to each
$(1, {1 \over k})$ with probability~$b_k$, where the $b_k$ are nonnegative
and sum to~1, and will be specified later. (Thus, the adversary's
actions are chosen to still be time-homogeneous Markov.)

To complete the construction, we choose $\{\beta_k\}$ and $\{a_k\}$
and $\{b_k\}$ so that
$\sum_k a_k r_k < \infty$ but $\sum_k b_k r_k = \infty$.
For example, we can do this by first choosing $\beta_k$
so that $r_k k^{-k} \to1$, and then letting
$a_k \propto(2k)^{-k}$ and $b_k \propto(k/2)^{-k}$.

It then follows that for the
Markov chain $P$ (governed by the $\{a_k\}$) the expected return time
to $O$ from $O$ is finite, and hence the chain
has a unique stationary probability measure $\pi$.
On the other hand, for the adversarial process $\{X_n\}$ (governed
by the $\{b_k\}$) the expected return time to $O$ from $O$ is
infinite. Hence, the adversarial process
is null recurrent, so it will move to larger
and larger $S_k$ as time progresses. In particular, the adversarial process
will \textit{not} be bounded in probability, even though the
transition probabilities $P$ are continuous.
\end{exa}

%re2 #&#
\begin{rem}
Example~\ref{exa2} is only defined on a countable state space
$\mathcal{X}$,
but if desired it could be ``extended'' to a counterexample on all of
$\mathbf{R}^2$.
For instance, we could let $\delta=10^{-6}$, and replace $\pi(\cdot
)$ by the
convolution $\pi(\cdot) * N(O,\delta^2)$ with a tiny normal distribution,
and replace $P(x,\cdot)$ by the
convolution $P(x,\cdot) * N(O,\delta^2)$ for each $x\in\mathcal{X}$,
and then continuously interpolate new
transition probabilities $P(x,\cdot)$ at all $x \in\mathbf{R}^2
\setminus\mathcal{X}$
such that $P(x,\cdot)$ is a probability measure for each $x\in\mathbf{R}^2$,
and the mapping $x \mapsto P(x,A)$ is continuous over $x\in\mathbf
{R}^2$ for
each fixed $A\in\mathcal{F}$.
This could be done in such a way that
(\ref{contcond}) would still be satisfied,
but (\ref{boundedinprob}) would still fail, thus providing
a counter-example even on the continuous state space $\mathbf{R}^2$.
\end{rem}

Finally, in a rather different direction, we consider what happens
if the process is allowed to be anticipatory, that is, to make moves
based on future randomness, with (\ref{XfollowsP}) replaced by the
weaker condition that $\P(X_{n+1} \in A | X_n=x) = P(x,A)$ but
without conditioning on the previous history $X_0,\break \ldots,X_{n-1}$.
% (Intuitively, the adversary is allowed to be ``clairvoyant''.)
It turns out that, under this subtle change, our theorems no longer hold:

%ex3 #&#
\begin{exa}\label{exa3}
Let $\mathcal{X}= [0,\infty) \subseteq\mathbf{R}$.
Define Markov chain transitions $P$ as follows.
For $x \le1$, $P(x,\cdot) =  \operatorname{Uniform}[0,2]$.
For $1 < x \le3$, $P(x,\cdot) = \break\operatorname{Uniform}[x-1,x+1]$.
For $3 < x \le4$, $P(x,\cdot) = \operatorname{Uniform}[4,5]$.
For $x > 4$, $P(x,\cdot) = \frac{1}{2} \delta_{x+1}(\cdot)
+ \frac{1}{2} \operatorname{Uniform}[x-2,x-1]$, where $\delta_{x+1}$ is
a point-mass at $x+1$.
Then $P$ is $\phi$-irreducible, with
negative drift for $x>4$, so $P$ must be
positive recurrent with some stationary probability
distribution $\pi$ to which it converges as in~(\ref{eqn-ergodic}).
Also, $P$~never moves more than a distance $D=2$
as in (\ref{eqn-boundedjumps}).

We next define the adversarial process $\{X_n\}$.
Let $K=[0,2]$, so $K_D=[0,4]$ and $K_{2D}=[0,6]$.
Let $\{B_i\}_{i=0}^\infty$ be i.i.d. with $\P(B_i=0)=\P(B_i=1)=1/2$,
and let $\{U_i\}_{i=0}^\infty$ be i.i.d. $\sim\operatorname{Uniform}[0,1]$,
and let\vadjust{\goodbreak} $a_* = 4 + \sum_{i=1}^\infty B_i 2^{-i}$.
For any $r\in\mathcal{X}$, let $r[i]$ be the coefficient of $2^i$
in the nonterminating binary expansion of $r$,
so that $r = \sum_{i\in\mathbf{Z}} r[i] 2^i$.
Conditionally on $X_n$, we construct $X_{n+1}$ by:
(a) if $X_n \le1$, then $X_{n+1} = 2 U_n$;
(b) if $1 < X_n \le3$, then $X_{n+1} = X_n - 1 + 2 U_n$;
(c) if $3 < X_n \le4$, then $X_{n+1} = a_*$;
(d) if $X_n > 4$, then $X_{n+1} = I_n(X_n+1) +
(1-I_n)(X_n-1-U_n)$, where
$I_n = \mathbf{1}_{X_n[-n] = B_n}$
is the indicator function of whether the coefficient
of $2^{-n}$ in the binary expansion of $X_n$ is equal to $B_n$.

Then it is easily checked that $\{X_n\}$ follows the one-step
transitions $P$ for all $x\in\mathcal{X}$ (including $x \in K$), in
the sense
that $\P(X_{n+1} \in A | X_n=x) = P(x,A)$ for all $A$ (but
without also conditioning on $X_0,\ldots,X_{n-1}$).
Furthermore,
(A1) holds with $M=1$ and $\mu_* = \operatorname{Uniform}[4,5]$.
Also, (A2) holds for the same $\mu_*$ due to
$P$'s negative drift for $x>4$.

On the other hand, by construction $a_*$ has the property that
$a_*[-n]=B_n$ for all $n\in\mathbf{N}$. Hence, once the chain hits
the interval
$(3,4]$, then it will move to $a_*$, and from there it will always
add~1 with probability~1. Therefore, $X_n \to\infty$ with
probability~1, so $\{X_n\}$ is not bounded in probability,
so (\ref{boundedinprob}) does not hold. This process thus provides
a counterexample to Theorem~\ref{thm-cemetery} if we assume only that
$\P(X_{n+1} \in A | X_n=x) = P(x,A)$, without also
conditioning on the previous history $X_0,\ldots,X_{n-1}$
as in (\ref{XfollowsP}).\vspace*{-2pt}
\end{exa}

%s5 #&#
\section{Proof of Theorem~\texorpdfstring{\protect\ref{thm-cemetery}}{4} and
Proposition~\texorpdfstring{\protect\ref{discreteprop}}{2}}
\label{sec-cemproof}

We begin by letting $\{Y_n\}$ be a ``cemetery process'' which begins in
the distribution $\mu_*$ at time~0, and then follows the
fixed transition kernel $P$, and then \textit{dies} as soon as it hits $K_D$.
Assumption (A2) then says that
this cemetery process $\{Y_n\}$ has finite expected lifetime.
For $L > \ell_0 := \sup\{\eta(x,\mathbf{0}) \dvtx x \in K_D\}$,
let $B_L = \{x\in\mathcal{X}\dvtx\eta(x,\mathbf{0}) \ge L\}$, and
let $N_L$ denote the cemetery process's
total occupation time of $B_L$ (i.e., the number of
iterations that $\{Y_n\}$ spends in $B_L$ before it dies).
We then\vspace*{-2pt} have:

%le8 #&#
\begin{lem}
\label{lemma-cembound}
Let $\{X_n\}$ be the adversarial process as defined
previously. Then assuming \textup{(A1)}, for any $n\in\mathbf{N}$, and
any $L > \ell_0$, and any $x \in K$, we\vspace*{-2pt} have
\[
\P(X_n \in B_L | X_0=x) \le M
\mathbf{E}(N_L),
\]
where $N_L$ is the occupation time of $B_L$ for
the cemetery process $\{Y_n\}$ defined above.\vspace*{-6pt}
\end{lem}

\begin{pf}
Let $\sigma$ be the last return time of $\{X_n\}$ to $K_D$ by time $n$
(which must exist since $X_0 \in K_D$),
and let $\mu_k$ be the (complicated) law of $X_k$ when starting
from $X_0=x_0$. Then letting $I = K_D \setminus K$ (``inside'')
and $O = K_{2D} \setminus K_D$ (``outside''), we have\vspace*{-4pt}
\begin{eqnarray*}
&& \P(X_n \in B_L | X_0=x_0)\\[-3pt]
 &&\qquad=
\sum_{k=0}^{n-1} \P(X_n \in
B_L, \sigma=k | X_0=x_0)
\\[-3pt]
&&\qquad= \sum_{k=0}^{n-1} \int
_{y\in I} \int_{z\in O} \P(X_k \in
dy, X_{k+1} \in dz, X_n \in B_L, \sigma=k |
X_0=x_0)
\\[-2pt]
&&\qquad= \sum_{k=0}^{n-1} \int
_{y\in I} \int_{z\in O} \mu_k(dy)
P(y,dz)\\
&&\hspace*{49pt}\qquad\qquad{}\times \P(X_n \in B_L, \sigma= k | X_0=x_0,
X_k=y, X_{k+1}=z)
\\[-2pt]
&&\qquad \le \sum_{k=0}^{n-1} \int
_{y\in I} \int_{z\in O} \mu_k(dy) M
\mu_*(dz)\\[-2pt]
&&\hspace*{49pt}\qquad\qquad{}\times \P(X_n \in B_L, \sigma= k |
X_0=x_0, X_k=y, X_{k+1}=z)
\\[-2pt]
&&\qquad\le \sum_{k=0}^{n-1} \int
_{y\in I} \int_{z\in O} \mu_k(dy) M
\mu_*(dz) \P(Y_{n-k-1} \in B_L | Y_0=z)
\end{eqnarray*}
(by letting $Y_n = X_{n+k+1}$, and noting that if $\sigma=k$, then the
process did not return to $K_D$ by time $n$, so it behaved like the
cemetery process between times $n-k-1$ and $n$)\vspace*{-8pt}
\begin{eqnarray*}
&\le& M \sum_{k=0}^{n-1} \int
_{z\in O} \P(Y_{n-k-1} \in B_L |
Y_0=z) \mu_*(dz)
\\[-2pt]
& \le& M \sum_{j=0}^\infty\int
_{z\in O} \P(Y_{j} \in B_L |
Y_0=z) \mu_*(dz).
\end{eqnarray*}
But this last sum is precisely the expected total
number of iterations that the cemetery process $\{Y_n\}$ spends in $B_L$
when started from the distribution $\mu_*$.
\end{pf}

\begin{pf*}{Proof of Theorem~\protect\ref{thm-cemetery}}
For each $A\in\mathcal{F}$,
let $\nu(A)$ be the above cemetery process's
expected occupation measure, that is, the
expected number of iterations that the cemetery process $\{Y_n\}$ spends
in the subset $A$. Then the total measure $\nu(\mathcal{X})$
equals the expected lifetime of the cemetery process, and is thus
finite by~(A2). Hence, by the usual continuity of measures,
\[
\lim_{L\to\infty} \nu(B_L) = \nu\biggl(\bigcap
_L B_L \biggr) = \nu(\varnothing) = 0.
\]
This shows that $\mathbf{E}(N_L) \to0$ as $L\to\infty$.
Hence, by Lemma~\ref{lemma-cembound},
\[
\lim_{L\to\infty} \sup_{n\in\mathbf{N}} \P(X_n \in B_L | X_0=x_0)
\le M \lim_{L\to\infty} \mathbf{E}(N_L) = 0,
\]
so $\{X_n\}$ is bounded in probability.
\end{pf*}

We now turn our attention to discrete chains as in
Proposition~\ref{discreteprop}. We begin with a lemma.
[Here and throughout, $\mathbf{E}_x(\cdots)$ means expected value
conditional on the process starting at the initial state
$x\in\mathcal{X}$.]

%le9 #&#
\begin{lem}
\label{discretelemma}
For an
irreducible Markov chain on a discrete state space with stationary
probability distribution $\pi$, for any two states $x$ and $y$, we have
$\mathbf{E}_x(\tau_y)<\infty$; that is, the chain will move from $x$
to $y$ in
finite expected time.
\end{lem}

\begin{pf}
If this were not the case, then it
would be possible from $y$ to travel to $x$ and then take infinite
expected time to return to $y$. This would imply
that $\mathbf{E}_y(\tau_y)=\infty$, contradicting the
fact that we must have $\mathbf{E}_y(\tau_y) = 1/\pi(y) < \infty$ by
positive recurrence.
\end{pf}

\begin{pf*}{Proof of Proposition~\protect\ref{discreteprop}}
Since $\mathcal{X}$ is countable and $P$ is irreducible, $\pi(x)>0$
for all
$x\in\mathcal{X}$.
Let $O = K_{2D} \setminus K_D$, and assume that $\pi(O)>0$
[otherwise increase $D$ to make this so, which can be done
unless $\pi(K_D)=1$ in which case the statement is trivial].

Since $K_{2D}$ is finite, assumption~(A1)
with $\mu_* = \pi|_{K_{2D} \setminus K_D}$
follows immediately with, for example, $M = (\max_{x,z\in K_{2D}}
P(x,z)) / (\min_{z \in K_{2D}} \pi(z)) < \infty$.

Next, note that
$\mathbf{E}_x(\tau_{K_D}) < \infty$ for each individual $x\in O$;
indeed, this follows by applying Lemma~\ref{discretelemma}
with any one specific $y\in K_D$
(which must exist since we assume $K$ is nonempty).
But then $\mathbf{E}_{\mu_*}(\tau_{K_D}) = \sum_{x \in O}
\mu_*(x) \mathbf{E}_x(\tau_{K_D})$, which must also be finite
since $O$ is finite. Hence, (A2) also holds.
The result thus follows from Theorem~\ref{thm-cemetery}.
\end{pf*}

%s6 #&#
\section{Two additional probability lemmas}

In this section, we prove two probability results which we
will use in the following section.

We first consider expected hitting times. Lemma~\ref{discretelemma}
above shows that \textit{discrete} ergodic Markov chains always
have $\mathbf{E}_x(\tau_y)<\infty$. On a general state space, one might
think by analogy that for any positive-recurrent $\phi$-irreducible
Markov chain with stationary distribution $\pi$, if $\pi(A)>0$ and
$\pi(B)>0$, then we must have $\mathbf{E}_{\pi|A}(\tau_B) < \infty
$, where
$\tau_B$ is the hitting time of $B$.
However, this is false. For example, consider a birth-death chain
on the positive integers having stationary distribution $\pi(j)
\propto j^{-2}$. Then if $B=\{1\}$ and $A=\{J,J+1,J+2,\ldots\,\}$ for
any $J > 1$, then $\mathbf{E}_{\pi|A}(\tau_B) \ge\sum_{j=J}^\infty
\pi(j)
(j-1) \propto\sum_{j=J}^\infty j^{-2} (j-1) = \infty$.

On the other hand, this result \textit{is} true in the case $A=B$.
Indeed, we have:

%le10 #&#
\begin{lem}
\label{lemma-expret}
Consider a Markov chain with stationary
probability distribution~$\pi$, and let $A\in\mathcal{F}$
with $\pi(A)>0$. Then:
\begin{longlist}[(ii)]
\item[(i)] $\mathbf{E}_{\pi|A}(\tau_A) = 1/\pi(A) < \infty$, where
$\tau_A$ is the first return time to $A$.
\item[(ii)]
For all $k\in\mathbf{N}$,
$\mathbf{E}_{\pi|A}(\tau_A^{(k)}) = k/\pi(A) < \infty$,
where $\tau_A^{(k)}$ is the $k$th return time to $A$.
\end{longlist}
\end{lem}

\begin{pf}
Part~(i) is essentially the formula of Kac~\cite{kac}. Indeed, using
Theorem~10.0.1 of \cite{MT} with $B=\mathcal{X}$, we obtain
\[
1 = \pi(\mathcal{X}) = \int_{x \in A} \pi(dx)
\mathbf{E}_x \Biggl[ \sum_{n=1}^{\tau_A}
\mathbf{1}_{X_n \in\mathcal{X}} \Biggr] = \int_{x \in A} \pi(dx)
\mathbf{E}_x[\tau_A] = \pi(A) \mathbf{E}_{\pi|A}[
\tau_A],
\]
giving the result.

For part~(ii), we
expand the original Markov chain to a new Markov chain on $\mathcal
{X}\times
\{0,1,\ldots,k-1\}$, where the first variable is the original chain,
and the second
variable is the count (mod $k$) of the number of times the chain
has returned to $A$. That is, each time the original chain visits $A$,
the second variable increases by 1 (mod $k$). Then the expanded chain
has stationary distribution $\pi\times\operatorname{Uniform}\{
0,1,\ldots,k-1\}$.
Hence, by part~(i), if we begin in
$(\pi|A) \times\delta_0$, then the expected return time of the
expanded chain to
$A \times\{0\}$ equals $1 / [\pi(A) \times(1/k)] = k/\pi(A)$.
But the first return time of the expanded chain to $A \times\{0\}$
corresponds precisely to the $k$th return time of the
original chain to $A$.
\end{pf}

We also require the following generalization of Wald's equation.

%le11 #&#
\begin{lem}
\label{waldgen}
Let $\{W_n\}$ be a sequence of nonnegative random variables
each with finite mean $m<\infty$,
and let $\{I_n\}$ be a sequence of indicator variables
each with $\P(I_n=1) = p > 0$.
Assume that the sequence of pairs  $\{(W_n,\break I_n)\}$ is i.i.d.
[i.e., the sequence\hspace*{1pt} $\{Z_n\}$ is i.i.d. where $Z_n = (W_n,I_n)$].
Let $\tau= \inf\{n \dvtx  I_n=1\}$, and let $S = \sum_{i=1}^\tau W_i$.
Then $\mathbf{E}(S) = {m \over p} < \infty$.
\end{lem}

\begin{pf}
We can write $S = \sum_{i=1}^\infty W_i \mathbf{1}_{\tau\ge i}$.
Now, the event $\{\tau\ge i\}$ is equivalent to the event that
$I_1=I_2=\cdots=I_{i-1}=0$.
Hence it is contained in $\sigma(Z_1,\ldots,Z_{i-1})$ and is thus
independent of $W_i$ by assumption.
Also, $\tau$ is distributed as Geometric($p$) and hence has mean $1/p$.
We then compute that
\begin{eqnarray*}
\mathbf{E}(S) &=& \mathbf{E}\Biggl( \sum_{i=1}^\infty
W_i \mathbf{1}_{\tau\ge i} \Biggr) = \sum
_{i=1}^\infty\mathbf{E}(W_i
\mathbf{1}_{\tau\ge i})
\\
&=& \sum_{i=1}^\infty\mathbf{E}(W_i)
\mathbf{E}(\mathbf{1}_{\tau
\ge i}) = \sum_{i=1}^\infty
m \P(\tau\ge i) = m \mathbf{E}(\tau) = m / p,
\end{eqnarray*}
as claimed.
\end{pf}

% \proof
% If $p=1$ then $\tau=1$ and the statement is trivial,
% so assume $p \in(0,1)$. We have
% $$
% m
% \ = \ \E[W_i]
% \ = \ p \E[W_i | I_i=1]
% + (1-p) \E[W_i | I_i=0]
% ,
% $$
% so $\E[W_i | I_i=1] \le m/p$ and
% $\E[W_i | I_i=0] \le m/(1-p)$.
% We then have that
% $$
% \E(S)
% \ = \ \sum_{n=1}^\infty
% \P[\tau=n] \ \E\Big[ \sum_{i=1}^n W_i \ \Big| \ \tau=n \Big]
% $$
% $$
% \ = \ \sum_{n=1}^\infty
% (1-p)^{n-1} p
% \ \ \E\Big[ \sum_{i=1}^n W_i \ \Big| \ I_1=\ldots=I_{n-1}=0, \ I_n=1
%\Big]
% $$
% $$
% \ = \ \sum_{n=1}^\infty
% (1-p)^{n-1} p \ \Big[ \Big( \sum_{i=1}^{n-1} \E[W_i | I_i=0] \Big)
% + \E[W_n | I_n=1] \Big]
% $$
% $$
% \ \le\ \sum_{n=1}^\infty
% (1-p)^{n-1} p \ \Big[ (n-1)m/(1-p) + m/p \Big]
% $$
% $$
% \ \le\ \Big[ m/\min(p,1-p) \Big]
% \ \sum_{n=1}^\infty n (1-p)^{n-1} p
% .
% $$
% This last sum is a finite expected value, thus giving the result.
% \qed

%s7 #&#
\section{Proof of Theorem~\texorpdfstring{\protect\ref{thm-pismall}}{5}}
\label{sec-smallproof}

The key to the proof is the following fact about Markov chain hitting
times.

%le12 #&#
\begin{lem}
\label{MChitting}
Consider a $\phi$-irreducible Markov chain on a state space
$(\mathcal{X},\mathcal{F})$ with transition
kernel $P$ and stationary probability distribution $\pi$. Let $B,C\in
\mathcal{F}$
with $\pi(B)>0$ and $\pi(C)>0$,
and let $\mu$ be any probability measure on $(\mathcal{X},\mathcal{F})$.
Suppose $C$ is a \textrm{small set} for $P$ with minorizing measure $\mu
$; that is, there is $\epsilon>0$ and $n_0\in\mathbf{N}$ such that
$P^{n_0}(x,A) \ge \epsilon \mu(A)$ for all
states $x \in C$ and all subsets $A \in\mathcal{F}$.
Let $\tau_B$ be the first hitting time of $B$. Then
$\mathbf{E}_\mu(\tau_B)<\infty$.
\end{lem}

\begin{pf}
It suffices to consider the case where $n_0=1$, since if not we
can replace $P$ by $P^{n_0}$ and note that the hitting time of $B$
by $P$ is at most $n_0$ times the hitting time of $B$ by~$P^{n_0}$.

We use the Nummelin splitting technique \cite{nummelin,MT}.
Specifically, we expand the state space
to $\mathcal{X}\times\{0,1\}$, where the second variable is an
indicator of
whether or not we are currently regenerating according to $\mu$.

Let $\alpha= \mathcal{X}\times\{1\}$. Then $\alpha$ is
a Markov chain atom (i.e., the chain has identical transition
probabilities from every state in $\alpha$), and it has stationary
measure $\pi(\alpha) = \epsilon \pi(C) > 0$. So,
by Lemma~\ref{lemma-expret}(i) above, if the expanded chain is
started in $\alpha$ (corresponding to the original chain starting
in $\mu$), then it will return to
$\alpha$ in finite expected time $1/\pi(\alpha) < \infty$.

We now let $W_n$ be the number of iterations between the
$(n-1)$st and $n$th returns to $\alpha$, and let
$I_n=1$ if this $n$th tour visits $B$, otherwise $I_n=0$.
Then $\P[I_n=1] > 0$ by the $\phi$-irreducibility of $P$.
Hence, $\{(W_n,I_n)\}$ satisfies the conditions of Lemma~\ref{waldgen}.

Therefore, by Lemma~\ref{waldgen}, the expected number of iterations
until we complete a tour which includes a visit to $B$ is finite.
Hence, the expected hitting time of $B$ is finite.
\end{pf}

%co13 #&#
\begin{cor}
\label{Aimply}
\textup{(A3)} with $\nu_*=\mu_*$ implies \textup{(A2)}.
\end{cor}

\begin{pf}
This follows immediately by applying Lemma~\ref{MChitting} with
$C = K_{2D} \setminus K_D$, and $B = K_D$, and $\mu= \mu_* = \nu_*$.
\end{pf}

\begin{pf*}{Proof of Theorem~\protect\ref{thm-pismall}}
Under assumption (a) that $\nu_*=\mu_*$, result (\ref{boundedinprob})
follows by combining Corollary~\ref{Aimply} with Theorem~\ref{thm-cemetery}.
Under assumption (b) that $P$ is reversible
and $\mu_*=\pi_{K_{2D} \setminus K_D}$,
it follows from the \hyperref[sec-appendix]{Appendix} (Section~\ref{sec-appendix} below)
that (A3) also holds with
$\nu_* = \pi|_{K_{2D} \setminus K_D} = \mu_*$.
Hence assumption (a) still applies, so (\ref{boundedinprob}) again follows.
\end{pf*}

%re3 #&#
\begin{rem}
One might wonder if it suffices in Theorem~\ref{thm-pismall} to
assume (A1) with \textit{any} distribution $\mu_*$,
and (A3) with \textit{any} distribution $\nu_*$,
without requiring that either $\nu_*=\mu_*$ or
$\mu_*=\pi|_{K_{2D} \setminus K_D}$.
Under these assumptions, it would still follow from
Lemma~\ref{lemma-expret}(ii) that the return times to
$K_{2D}$ all have finite expectation. And it would still be true
that \textit{if} we regenerate from $\nu_*$ in finite expected
time, then we will eventually hit $K_D$ in finite expected time.
The problem is that the expected
time to \textit{first} regenerate from $\nu_*$ might be infinite.
Indeed, conditionally upon visiting $K_{2D}$ but repeatedly
\textit{failing} to regenerate, the chain could perhaps move to
worse and worse states from which it would then take longer and
longer to return to $K_{2D}$. (It~is tempting to apply
Lemma~\ref{waldgen} here where $W_n$ is the time between consecutive
visits to $K_{2D}$ and $I_n=1$ if we regenerate otherwise 0, but
unfortunately in this case $\{(W_n,I_n)\}$ are not i.i.d., and conditionally
on nonregeneration the values of $\mathbf{E}[W_n | I_1=\cdots=I_n=0]$
could grow unboundedly.)
\end{rem}

%s8 #&#
\section{Proof of Proposition~\texorpdfstring{\protect\ref{unifsmallprop}}{6}}
\label{sec-unifsmallproof}

Let $A \subseteq\mathbf{R}^d$ be the ball centered at the origin of
radius~1, and let $B \subseteq\mathbf{R}^d$ be the ball centered at the
point $(3/2,0,0,\ldots,0)$ of~radius~1. Then $A \cap B$ has nonempty
interior, so $v_d := \operatorname{Leb}(A \cap B) > 0$. In terms of
this, we have:

%le14 #&#
\begin{lem}\label{balloverlap}
Let $A,B \subseteq\mathbf{R}^d$ be two balls with radii $r \le
R$, such that their centers are a distance $w \le3r/2+(R-r)$ apart. Then
$\operatorname{Leb}(A \cap B) \ge r^d v_d$.
\end{lem}

\begin{pf}
If $r=R=1$, then this is just the definition of $v_d$.
If one of the balls is stretched by a factor $R>1$ while moving its
center a distance $R-r$ further away, then the new ball contains
the old ball, so $\operatorname{Leb}(A \cap B)$ can only increase.
Finally, if each of $w$ and $r$ and $R$ are multiplied by the same
constant $a>0$, then the entire geometry is scaled by a factor of $a$,
so $\operatorname{Leb}(A \cap B)$ is multiplied by $a^d$.
Combining these facts, the result follows.
\end{pf}

%le15 #&#
\begin{lem} \label{Plower}
Let $P$ be a Markov chain on an open subset $\mathcal{X}\subseteq
\mathbf{R}^d$.
Let $J$ be
a rectangular subset of $\mathcal{X}$, of the form $J = (a_1,b_1)
\times\cdots
\times(a_d,b_d) \subseteq\mathcal{X}$, where
$a_i<b_i$ are extended real numbers
(i.e., we might have $a_i=-\infty$ and/or
$b_i=\infty$ for some of the $i$).
Suppose there are $\delta>0$ and $\epsilon>0$ satisfying the
condition~(\ref{deltepeqn}) that
$P(x,dy) \ge\epsilon \operatorname{Leb}(dy)$ whenever $x,y \in J$
with $|y-x|<\delta$.
Then for each $n\in\mathbf{N}$, there is $\beta_n>0$ such that
$P^n(x,dy) \ge\beta_n \operatorname{Leb}(dy)$ whenever $x,y \in J$
with $|y-x| < \delta(n+1)/2$.
\end{lem}

\begin{pf}
We first consider the case where $a_i=-\infty$ and $b_i=\infty$ for
all~$i$.
The result for $n=1$ follows by assumption.
Suppose the result is true for some $n \ge1$.
Let $|y-x|<\delta(n+1)/2$,
let $A$ be the ball centered at $x$ of radius $\delta(n+1)/2$
and let $B$ be the ball centered at $y$ of radius $\delta$.
Then applying Lemma~\ref{balloverlap} with $r=\delta$
and $R=\delta(n+1)/2$ and $w=\delta(n+2)/2$, we see that
$\operatorname{Leb}(A \cap B) \ge\delta^d v_d$.
The result now follows from the calculation
\begin{eqnarray*}
P^{n+1}(x,dy) &=& \int_{z\in\mathcal{X}} P^n(x,dz)
P(z,y) \ge \int_{z\in A \cap B} P^n(x,dz) P(z,y)
\\
&\ge& \int_{z\in A \cap B} \beta_n \operatorname{Leb}(dz)
\epsilon \operatorname{Leb}(dy) \ge \operatorname{Leb}(A \cap B)
\beta_n \epsilon \operatorname{Leb}(dy)
\\
&\ge& \delta^d v_d \beta_n \epsilon
\operatorname{Leb}(dy) =: \beta_{n+1} \operatorname{Leb}(dy).
\end{eqnarray*}

For the general case, by shrinking $\delta$ as necessary, we can assume
that $\delta< \frac{1}{2} \min_i (b_i-a_i)$. Then in the above calculation
we can only use those parts of $A \cap B$ which
are still inside $J$. But here $J$ must contain at least half of $A
\cap B$ in each coordinate, that is, at least $1/2^d$ of $A \cap B$ overall.
Hence, $\operatorname{Leb}(A \cap B \cap J) \ge(1/2^d) \operatorname
{Leb}(A \cap B)$.
So, the above calculation still goes through, except now
with $\beta_{n+1}$ multiplied by an extra factor of $1/2^d$.
\end{pf}

\begin{pf*}{Proof of Proposition~\protect\ref{unifsmallprop}}
Let $z = \operatorname{Diam}(J) < \infty$.
Find $n_0\in\mathbf{N}$ such that $\delta(n_0+1)/2 > z$.
Then it follows from Lemma~\ref{Plower} that
there is $\epsilon_{n_0}>0$ such that
$P^{n_0}(x,dy) \ge\epsilon_{n_0} \operatorname{Leb}(dy)$ for all
$x,y \in J
\supseteq K_{2D} \setminus K_D$.
Hence, (A3) holds for this $n_0$ with
$\nu_* = \operatorname{Uniform}(K_{2D} \setminus K_D)$ and
$\epsilon= \epsilon_{n_0} \operatorname{Leb}(K_{2D} \setminus K_D)$.
\end{pf*}

%s9 #&#
\section{Some facts about geometric ergodicity}

To prove Theorem~\ref{thm-geomerg}, we need to understand the
implications of the geometric ergodicity assumption (A4).
The following proposition shows that we can always find a geometric
drift function of a certain form.
To state it, let $PV(x) = \int_{y\in\mathcal{X}} V(y) P(x,dy)$ be
the action
of the Markov kernel $P$ on a function $V$, and let
$\tau_C = \inf\{n \ge1 \dvtx X_n \in C\}$ be the first hitting time
of $C$
by a Markov chain $\{X_n\}$ following the transitions $P$.
Also, say that $V$ is a \textit{geometric drift function} if
%
%e7 #&#
\begin{equation}
\label{drifteqn} PV(x) \le \lambda V(x) + b \mathbf{1}_C(x)
\end{equation}
for some small set $C\in\mathcal{F}$ and some real numbers
$\lambda<1$ and $b<\infty$.

%pr16 #&#
\begin{prop}
\label{geomfacts}
If $P$ is geometrically ergodic as in~(A4), then there is a
small set $C \subseteq\mathcal{X}$ with $\pi(C)>0$,
and a real number $\kappa>1$, such that
the function $V\dvtx\mathcal{X}\to\mathbf{R}$ defined by $V(x) =
\mathbf{E}_x(\kappa^{\tau_C})$
is $\pi$-a.e. finite, and $r := \sup_{x \in C} V(x) < \infty$,
and the geometric drift equation (\ref{drifteqn}) holds
with this $C$ for some $b<\infty$
and with $\lambda= \kappa^{-1} < 1$.
Furthermore, there is $\rho<1$ and $c<\infty$ such that
$\|P^n(x,\cdot) - \pi\| \le c V(x) \rho^n$
for all $n\in\mathbf{N}$ and $x\in\mathcal{X}$.
\end{prop}

\begin{pf}
Let $A_M = \{x\in\mathcal{X}\dvtx\xi(x) \le M\}$.
Since $\pi\{x\in\mathcal{X}\dvtx\xi(x)<\infty\}=1$, we can find
$M<\infty$ with $\pi(A_M)>0$.
The existence of \textit{some} small set $C \subseteq A_M$ with $\pi(C)>0$
follows from, for example, \cite{orey} (where they are called \textit
{$C$-sets})
or \cite{nummelin}
or Theorem~5.2.2 of \cite{MT}.
The fact that $C \subseteq A_M$ then implies condition~(15.1) of \cite{MT}
for this $C$ [with $P^\infty(C)=\pi(C)$ and
$M_C=M$ and $\rho_C=\rho$]. The existence of a
(possibly different) small set $C$ and $\kappa>1$
with $\pi(C)>0$ and $r := \sup_{x \in C} \mathbf{E}_x(\kappa^{\tau
_C}) < \infty$.
then follows from Theorem 15.0.1(ii) of \cite{MT}.

Let $V(x) = \mathbf{E}_x(\kappa^{\tau_C})$.
We compute directly that
if $\{W_n\}$ follows $P$, then for $x \notin C$,
\begin{eqnarray*}
V(x) &=& \mathbf{E}\bigl(\kappa^{\tau_C} | W_0=x\bigr) =
\mathbf{E}\bigl[ \mathbf{E}\bigl(\kappa^{\tau_C} | W_1\bigr) |
W_0=x \bigr]
\\
&=& \int_{y\in\mathcal{X}} \mathbf{E}\bigl(\kappa^{\tau_C} |
W_1=y\bigr) P(x,dy) = \int_{y\in\mathcal{X}} \mathbf{E}
\bigl(\kappa^{\tau_C+1} | W_0=y\bigr) P(x,dy)
\\
&=& \int_{y\in\mathcal{X}} \kappa \mathbf{E}\bigl(\kappa^{\tau_C}
| W_0=y\bigr) P(x,dy) = \kappa \int_{y\in\mathcal{X}} V(y)
P(x,dy) = \kappa PV(x),
\end{eqnarray*}
which shows that $PV(x) = \kappa^{-1} V(x)$ for $x \notin C$.

To prove the geometric drift condition,
it remains only to prove that $b := \sup_{x \in C} PV(x)$ is finite.
For this we use some additional results from \cite{MT}.
We first compute that in the special case $f \equiv1$, we\vspace*{-3pt} have that
\begin{eqnarray*}
\sup_{x \in C} \mathbf{E}_x \Biggl( \sum
_{k=0}^{\tau_C-1} f(W_k) \kappa ^k
\Biggr) &=& \sup_{x \in C} \mathbf{E}_x \Biggl( \sum
_{k=0}^{\tau_C-1} \kappa^k \Biggr)= \sup_{x \in C} \mathbf{E}_x \biggl(
{ \kappa^{\tau_C} - 1 \over
\kappa-1 } \biggr)
\\
 &=& { \sup_{x \in C} \mathbf{E}_x(\kappa^{\tau_C}) - 1 \over\kappa
-1 } = { r - 1 \over\kappa-1 }
< \infty.
\end{eqnarray*}
This means\vspace*{-1pt} that $C$ is an ``$f$-Kendall set'' for $f \equiv1$,
as defined on page~368 of~\cite{MT}. Hence,
by Theorem~15.2.4 of \cite{MT}, the function
$G(x) := G_C^{(\kappa)}(x,f)$ which equals~1 inside $C$ and\vspace*{-3pt} equals
%
%e8 #&#
\begin{equation}
\label{Geqn} \mathbf{E}_x \Biggl( \sum
_{k=0}^{\tau_C} f(W_k) \kappa^k
\Biggr) = \mathbf{E}_x \Biggl( \sum_{k=0}^{\tau_C}
\kappa^k \Biggr) = { \mathbf{E}_x(\kappa^{\tau_C+1}) - 1 \over\kappa- 1 } = { \kappa V(x) - 1 \over\kappa- 1 }
\end{equation}
outside of $C$, satisfies its own geometric drift condition,
say $PG(x) \le\break \lambda_G G(x) + b_G \mathbf{1}_C(x)$ where $\lambda_G<1$
and $b_G<\infty$.
In particular, since $G(x)=1$ for $x \in C$,
this means that $\sup_{x \in C} PG(x) \le\lambda_G + b_G < \infty$.
Now, by~(\ref{Geqn}), for $x\notin C$
we have $V(x) = {1 \over\kappa} [ 1 + (\kappa-1) G(x) ]
\le1+G(x)$.
Since $V(x) \le r$ for $x \in C$, it follows that
for all $x\in\mathcal{X}$, we have $V(x) \le r+G(x)$.
Therefore, $PV(x) \le r+PG(x)$. This shows, finally, that
\[
b := \sup_{x \in C} PV(x) \le r + \sup_{x \in C}
PG(x) \le r + \lambda_G + b_G < \infty.
\]

The above two facts together show that
$PV(x) \le\kappa^{-1} V(x) + b \mathbf{1}_C(x)$ with \mbox{$b<\infty$}.

The bound\vspace*{-3pt} on $\|P^n(x,\cdot) - \pi\|$
then follows from Theorem~16.0.1 of \cite{MT}.
\end{pf}

We next establish some bounds based on
geometric-drift-type\vspace*{-3pt} inequalities.

%le17 #&#
\begin{lem}
\label{gendriftlemma}
Let $\{Z_n\}$ be any stochastic process. Suppose there
are $0<\lambda<1$ and $b<\infty$ such that for all $n\in\mathbf
{N}$, we have
$
\mathbf{E}(Z_n | Z_0,\ldots,Z_{n-1})
\le \lambda Z_{n-1} + b
$.
Then\vspace*{-6pt} for all $n\in\mathbf{N}$,
\[
\mathbf{E}(Z_n | Z_0) \le \lambda^n
Z_0 + {b \over1-\lambda} \le Z_0 +
{b \over1-\lambda}.
\]
\end{lem}

\begin{pf}
We claim that for all $n \ge0$,
%
%e9 #&#
\begin{equation}
\label{inductionstep} \mathbf{E}(Z_n | Z_0) \le
\lambda^n Z_0 + \bigl(1+\lambda+\cdots+
\lambda^{n-1}\bigr) b.
\end{equation}
Indeed, for $n=0$ this is trivial, and for $n=1$ this
is equivalent to the hypothesis of the lemma.
Suppose now that (\ref{inductionstep}) holds for some value of $n$. Then
%
%e10 #&#
%e11 #&#
\begin{eqnarray*}
\mathbf{E}(Z_{n+1} | Z_0) &=& \mathbf{E}\bigl(
\mathbf{E}(Z_{n+1} | Z_0,\ldots,Z_n) |
Z_0 \bigr) \le \mathbf{E}( \lambda Z_n + b |
Z_0 )
\\[-2pt]
&\le& \lambda \bigl( \lambda^n Z_0 + \bigl(1+\lambda+
\cdots+\lambda^{n-1}\bigr) b \bigr) + b \\[-2pt]
&=& \lambda^{n+1}
Z_0 + \bigl(1+\lambda+\cdots+\lambda^{n-1}+
\lambda^n\bigr) b,
\end{eqnarray*}
so (\ref{inductionstep}) holds for $n+1$. Hence, by induction,
(\ref{inductionstep}) holds for all $n \ge0$.

The result now follows since
$1+\lambda+\cdots+\lambda^{n-1} = {1 - \lambda^n \over1 - \lambda}
\le{1 \over1 - \lambda}$.
\end{pf}

%pr18 #&#
\begin{prop}
\label{finiteexpprop}
If $P$ is geometrically ergodic with stationary probability
distribution $\pi$ and $\pi$-a.e. finite geometric drift function
$V$ satisfying $PV(x) \le\lambda V(x) + b$ where
$0 \le\lambda< 1$ and $0 \le b < \infty$, then
$\mathbf{E}_\pi(V) \le b/(1-\lambda) < \infty$.
\end{prop}

\begin{pf}
Choose any $x\in\mathcal{X}$ with $V(x)<\infty$ (which holds for
$\pi$-a.e. $x\in\mathcal{X}$).
Then applying Lemma~\ref{gendriftlemma} to $Z_n = P^nV(x)$ gives
$P^nV(x) \le V(x) + {b \over1-\lambda}$, and in particular
$P^nV(x) \not\to\infty$.
But Theorem~14.3.3 of \cite{MT} with $f=V$ states that if
$\pi(V)=\infty$, then $P^nV(x) \to\infty$ for all $x\in\mathcal{X}$.
Hence, by contraposition, we must have $\pi(V)<\infty$.

Finally, we have by stationarity that $\pi(V) = \pi(PV)$. So, taking
expectations with respect to $\pi$ of both sides of the inequality
$PV \le\lambda V + b$ and using that $\pi(V)<\infty$, we obtain that
$\pi(V) \le\lambda \pi(V) + b$, whence $\pi(V) \le b/(1-\lambda)$.
\end{pf}

%re4 #&#
\begin{rem}
If $P$ is \textit{uniformly
ergodic}, meaning that (A4) holds for a \textit{constant}
function $V<\infty$, then it follows from Theorem~16.0.2(vi) of
\cite{MT} that $U := \sup_{x\in\mathcal{X}} \mathbf{E}_x(\tau
_{K_D}) < \infty$,
which implies that $\mathbf{E}_{\mu_*}(\tau_{K_D}) \le U < \infty$,
so (A2) must hold.
\end{rem}

%s10 #&#
\section{Proof of Theorem~\texorpdfstring{\protect\ref{thm-geomerg}}{7}}
\label{sec-geomproof}

The key to the proof is a uniform bound on certain powers of $P$:

%le19 #&#
\begin{lem}
\label{supPVlemma}
Assuming \textup{(A4)} and \textup{(A5)}, with
$V$ as in Proposition~\ref{geomfacts},
$\sup_{x\in K_D} \sup_{n \ge0} P^nV(x) < \infty$.
\end{lem}

\begin{pf}
For $x\in K_D$, $PV(x) = \mathbf{E}_{y\sim P(x,\cdot)} V(y) \le M
\mathbf{E}_{y\sim\pi}V(y) = M \pi(V) < \infty$ by
(A5) and Proposition~\ref{finiteexpprop}.
Then applying Lemma~\ref{gendriftlemma} to $Z_n = P^nV(x)$ gives
$P^nV(x) \le M \pi(V) + {b \over1-\lambda}$.
In particular, $\sup_{x\in K_D} \sup_{n \ge1} P^nV(x) < \infty$.

Furthermore, for $x\in K_D$, $V(x) = \mathbf{E}_x(\kappa^{\tau_C})
= \kappa \mathbf{E}_{P(x,\cdot)}(\kappa^{\tau_C})
\le\kappa M \mathbf{E}_\pi(\kappa^{\tau_C})
= \kappa M \pi(V) < \infty$ by Proposition~\ref{finiteexpprop},
so the above ``sup'' can be extended to include $n = 0$ too.
\end{pf}

%re5 #&#
\begin{rem}
For Metropolis algorithms on continuous state spaces,
usually $P(x,\{x\})>0$ for most $x\in\mathcal{X}$, so
(A5) usually won't hold [though
(A1) often will; see Section~\ref{sec-mcmc}].
On the other hand, if
$
P(x,\cdot) = r(x) \delta_x(\cdot) + (1-r(x)) R(x,\cdot)
$
where $\delta_x$ is a point-mass at $x$ and $0 \le r(x) \le1$ and
$R$ satisfies (A5), then it is easily seen that if
$\kappa r(x) \le B < 1$ for all $x \in K_D$, then
Lemma~\ref{supPVlemma} still holds with
$\sup_{x \in K_D} V(x) \le \kappa M \pi(V) / (1 - B) < \infty$,
and the rest of the
proof of Theorem~\ref{thm-geomerg} then goes through without change.
\end{rem}

%pr20 #&#
\begin{prop}
\label{geomVprop}
Assuming \textup{(A4)} and \textup{(A5)}, the random sequence\break
$\{V(X_n)\}$ is bounded in probability.
\end{prop}

\begin{pf}
Lemma~\ref{supPVlemma} with $n=0$ says that
$U := \sup_{x\in K_D} V(x) < \infty$. Since
the adversary can only adjust the values of $\{X_n\}$ within $K_D$, it
follows that the adversary can only change the ``next value of
$V(X_n)$'' by at most $U$, so $\{X_n\}$ will still satisfy a drift
condition similar to (\ref{drifteqn}),
for the same $C$ and $\lambda$
but with $b$ replaced by $b+U < \infty$.
(Of course, $C$ might not be a small set for the adversarial process.)
More precisely, it follows from
(\ref{drifteqn}) that the adversarial process $\{X_n\}$ satisfies that
$\mathbf{E}[V(X_{n}) | X_0,X_1,\ldots,X_{n-1}] \le\lambda
V(X_{n-1}) + b + U$.
Hence, applying Lemma~\ref{gendriftlemma} to $Z_n = V(X_n)$ says that
$\{\mathbf{E}_{x_0}[V(X_n)]\}$ is bounded in probability,
that is, that\vspace*{1pt} $\zeta:= \sup_{x\in K_D} \sup_{n \ge0} \mathbf
{E}[V(X_n) | X_0=x_0]
< \infty$.
It then follows by Markov's inequality that
$\P_{x_0}[V(X_n) \ge R] \le\zeta/R$ for all $n$ and all $R>0$. Hence,
$\{V(X_n)\}$ is bounded in probability.
\end{pf}

%re6 #&#
\begin{rem}
$\!\!$Proposition~\ref{geomVprop} immediately implies a bound on the
\textit{$\epsilon$-\break convergence times} \cite{adapt.tex} defined by
$M_\epsilon(x) =
\inf\{n \ge1 \dvtx\|P^n(x,\cdot) - \pi(\cdot)\| \le\epsilon\}$.
Indeed, by Proposition~\ref{geomfacts} we have
$\|P^n(x,\cdot) - \pi\| \le c V(x) \rho^n$, whence
$M_\epsilon(x) \le\lceil\log(c V(x)/\epsilon)/\log(1/\rho)
\rceil$.
Since $\{V(X_n)\}$ is bounded in probability by
Proposition~\ref{geomVprop},
it follows that $\{M_\epsilon(X_n)\}$ is bounded in probability too;
see also the containment condition (\ref{containmenteqn}) below.
\end{rem}

\begin{pf*}{Proof of Theorem~\protect\ref{thm-geomerg}}
The bounded-jumps condition~(\ref{eqn-boundedjumps})
implies that the small set $C$
must be bounded (in fact, of diameter $\le2Dn_0$).
Let $r = \sup\{|x| \dvtx x \in C\} < \infty$. Then if $|x|>r$,
it takes at least $(|x|-r)/D$ steps to return to $C$ from $x$.
Hence, $V(x) \ge\kappa^{(|x|-r)/D}$. Therefore,
$|x| \le r + D \log(V(x))/\break \log(\kappa)$, so
$|X_n| \le r + D \log(V(X_n))/\log(\kappa)$.
But $\{V(X_n)\}$ is bounded in probability by Proposition~\ref{geomVprop}.
Hence, so is $\{X_n\}$.
\end{pf*}

%s11 #&#
\section{Application to adaptive MCMC algorithms}
\label{sec-mcmc}
\textit{Markov chain Monte Carlo (MCMC) algorithms} proceed by running
a Markov chain $\{X_n\}$ with stationary probability distribution
$\pi$, in the hopes that $\{X_n\}$ converges in total variation
distance to $\pi$, that is, that
%
%e12 #&#
\begin{equation}
\label{MCMCconveqn} \lim_{n\to\infty} \sup_{A\in\mathcal{F}} \bigl|
\P(X_n \in A) - \pi (A) \bigr| = 0,\qquad x \in\mathcal{X}, A \in
\mathcal{F}.
\end{equation}
If so then for large $n$,
the value of $X_n$ is approximately a ``sample''
from $\pi$. Such algorithms are hugely popular in, for example, Bayesian
statistical inference; for an overview see, for example,~\cite{galinbook}.

\textit{Adaptive MCMC algorithms} \cite{haario} attempt to speed
up the convergence~(\ref{MCMCconveqn}) and thus make MCMC more
efficient by modifying the Markov chain transitions during the run
(i.e., ``on the fly'') in a search for a more optimal chain; for a
brief introduction, see, for example,~\cite{galinart}. Such algorithms often
appear to work very well in practice (e.g., \cite{adaptex.tex,kohn,craiu,craiu2}). However, they are no longer Markov
chains (since
the adaptions typically depend on the process's entire history),
making it extremely difficult to establish mathematically that
the convergence~(\ref{MCMCconveqn}) will even be preserved (much
less improved). As a result, many papers either make the artificial
assumption that the state space $\mathcal{X}$ is compact (e.g., \cite
{haario,craiu,craiu2}), or prove the convergence~(\ref{MCMCconveqn})
using complicated mathematical arguments requiring strong and/or
uncheckable assumptions (e.g., \cite{atchade,AndrieuMoulines,adapt.tex,kohn,AtchadeFort,vihola,adapGibbs}),
or do not prove
(\ref{MCMCconveqn}) at all and simply hope for the best. It is
difficult to find simple, easily checked conditions which provably
guarantee the convergence (\ref{MCMCconveqn}) for adaptive MCMC
algorithms.

One step in this direction is in \cite{adapt.tex}, where it is proved
that convergence (\ref{MCMCconveqn})
is implied by two conditions. The first condition is
\textit{diminishing adaptation}, which says that the process adapts
less and less as time goes on; see (\ref{dimadeqn}) below.
The second condition is \textit{containment},
which says that the process's convergence times are bounded in
probability; see (\ref{containmenteqn}) below.
The first of these two conditions is usually easy to satisfy directly
by wisely designing the algorithm, so it is not of great
concern. However, the second condition is notoriously difficult
to verify (see, e.g., \cite{bai}) and thus a severe limitation
(though an essential condition; cf. \cite{adapFail}).
On the other hand, the containment
condition (\ref{containmenteqn}) is reminiscent of the boundedness in
probability property (\ref{boundedinprob}), which is implied
by our various theorems above. This suggests that our theorems
might be useful in establishing the containment condition
(\ref{containmenteqn}) for certain adaptive MCMC algorithms, as we
now explore.

%s11.1 #&#
\subsection{The adaptive MCMC setup}
\label{sec-mcmcsetup}

We define an adaptive MCMC algorithm
within the context of Section~\ref{sec-setup} as follows.
Let $\mathcal{X}$ be an open subset of $\mathbf{R}^d$ for some $d\in
\mathbf{N}$,
on which $\pi$ is some probability distribution.
Assume that for some compact index set $\mathcal{Y}$,
there is a collection $\{P_\gamma\}_{\gamma\in\mathcal{Y}}$
of Markov kernels on $\mathcal{X}$, each of which leaves $\pi$
stationary and
in fact is Harris-ergodic to $\pi$ as in~(\ref{eqn-ergodic}).
The adversary proceeds by choosing, at
each iteration $n$, an index $\Gamma_n \in\mathcal{Y}$ (possibly
depending on
$n$ and/or the process's entire history, though not on the future).
The process $\{X_n\}$ then moves at time $n$ according to the
transition kernel $P_{\Gamma_n}$, that is,
\[
\P(X_{n+1} \in A | X_n=x, \Gamma_n=\gamma,
X_0,\ldots,X_{n-1}, \Gamma_0,\ldots,
\Gamma_{n-1}) = P_\gamma(x,A).
\]
To reflect the bounded jump condition (\ref{eqn-boundedjumps}), we assume
there is $D<\infty$ with
%
%e13 #&#
\begin{equation}
\label{boundedPgamma} P_\gamma\bigl( x, \{y\in\mathcal{X}\dvtx|y-x| \le D\}
\bigr) = 1,\qquad x \in\mathcal{X}, \gamma\in\mathcal{Y}.
\end{equation}
To reflect that the adversary can only adapt inside $K$, we assume
that the $P_\gamma$ kernels are all equal outside of $K$, that is, that
%
%e14 #&#
\begin{equation}
\label{fixedPeqn} P_\gamma(x,A) = P(x,A),\qquad A\in\mathcal{F}, x \in
\mathcal{X}\setminus K,
\end{equation}
for some fixed Markov chain kernel $P(x,dy)$
also satisfying~(\ref{eqn-ergodic}).
We further assume that
%
%e15 #&#
\begin{equation}
\label{Pboundedabove} \qquad\exists M<\infty\quad \mbox{s.t.}\quad P(x,dy) \le M
\operatorname{Leb}(dy),\qquad x \in K_D \setminus K, z \in
K_{2D} \setminus K_D.\hspace*{-3pt}
\end{equation}
We also assume the $\epsilon$--$\delta$ condition~(\ref{deltepeqn}) that
$P(x,dy) \ge\epsilon \operatorname{Leb}(dy)$ whenever $x,y \in J$
with $|y-x|<\delta$, for some bounded rectangle $J$ with
$K_{2D} \setminus K_D \subseteq J \subseteq\mathcal{X}$.

We shall particularly focus on the case where each $P_\gamma$
is a Metropolis--Hastings algorithm. This means that $P_\gamma$
proceeds, given $X_n$, by first choosing a proposal state
$Y_{n+1} \sim Q_\gamma(X_n,\cdot)$ for some proposal kernel
$Q_\gamma(x,\cdot)$ having a density $q_\gamma(x,y)$ with respect to
$\operatorname{Leb}$. Then,\vspace*{1pt} with probability
$\alpha_\gamma(X_n, Y_{n+1})
:= \min[ 1, { \pi(Y_{n+1}) q_\gamma(Y_{n+1}, X_n)
\over\pi(X_n) q_\gamma(X_n, Y_{n+1}) } ]$
it \textit{accepts} this
proposal by setting $X_{n+1} = Y_{n+1}$. Otherwise, with probability
$1 - \alpha_\gamma(X_n, Y_{n+1})$, it \textit{rejects} this proposal
by setting $X_{n+1} = X_n$.
% (If $q_\gamma(Y_{n+1}, \ X_n) = q_\gamma(X_n, \ Y_{n+1})$
% then these two factors cancel.)
That is,
\[
P_\gamma(x,A) = r(x) \delta_x(A) + \int
_{y \in A} Q_\gamma(x,dy) \alpha_\gamma(x,y),
\]
where $\delta_x(\cdot)$ is a point-mass at $x$, and
$r(x) = 1 - \int_{y \in\mathcal{X}} Q_\gamma(x,dy) \alpha_\gamma(x,y)$
is the overall probability of rejecting.
Note that (\ref{boundedPgamma}), (\ref{fixedPeqn})
and (\ref{Pboundedabove}) are each automatically
satisfied for $P_\gamma$ and $P$ if the corresponding
equations are satisfied for corresponding $Q_\gamma$ and $Q$.

%s11.2 #&#
\subsection{An adaptive MCMC theorem}

Our theorem shall follow up on the result from \cite{adapt.tex}
that convergence (\ref{MCMCconveqn})
is implied by the twin properties of diminishing adaptation and containment.
\textit{Diminishing adaptation} says that the process adapts
less and less as time goes on, or more formally that
%
%e16 #&#
\begin{equation}
\label{dimadeqn} \lim_{n\to\infty} \sup_{x\in\mathcal{X}} \sup
_{A\in\mathcal{F}} \bigl|P_{\Gamma_{n+1}}(x,A) - P_{\Gamma_n}(x,A)\bigr| = 0
\qquad \mbox{in probability}.
\end{equation}
\textit{Containment} says that the process's convergence times are
bounded in probability, or more formally that
%
%e17 #&#
\begin{equation}
\label{containmenteqn} \bigl\{M_\epsilon(X_n,\Gamma_n)
\bigr\}_{n=1}^\infty \qquad\mbox{is bounded in probability},
\end{equation}
where $M_\epsilon(x,\gamma) = \inf\{n \ge
1 \dvtx\|P_\gamma^n(x,\cdot) - \pi(\cdot)\| \le\epsilon\}$ is the
\textit{$\epsilon$-convergence time}.
The containment condition (unlike diminishing adaptation) is
notoriously difficult to establish in practice
(see, e.g., \cite{bai}), but the theorems herein can help.
To state a clean theorem, we assume continuous densities, as follows:
\begin{longlist}[(A6)]
\item[(A6)] $\pi$ has a
continuous positive density function (with respect to $\operatorname{Leb}$),
and the transition probabilities $P_\gamma(x,dy)$ either (i)
have densities which are continuous
functions of $x$ and $y$ and $\gamma$,
or (ii) are Metropolis--Hastings algorithms whose proposal kernel
densities $q_\gamma(x,dy)$ are continuous functions
of $x$, $y$ and $\gamma$.
\end{longlist}

In terms of the above setup, we have:

%th21 #&#
\begin{thm}
\label{mcmcthm}
Consider an adaptive MCMC algorithm as in Section~\ref{sec-mcmcsetup},
on an open subset $\mathcal{X}$ of $\mathbf{R}^d$,
such that the kernels $P_\gamma$ (or the proposal kernels $Q_\gamma$
in the case of adaptive Metropolis--Hastings) have
bounded jumps as in (\ref{boundedPgamma}),
and no adaption outside of $K$ as in (\ref{fixedPeqn}),
with the fixed kernel $P$ (or a corresponding fixed proposal kernel $Q$)
bounded above as in (\ref{Pboundedabove}). We further assume
the $\epsilon$--$\delta$ condition~(\ref{deltepeqn}) for $P$, and
the continuous densities condition \textup{(A6)}.
Then the algorithm satisfies the containment condition~(\ref{containmenteqn}).
Hence, assuming diminishing adaptation~(\ref{dimadeqn}),
the algorithm converges in distribution to $\pi$
as in~(\ref{MCMCconveqn}).
\end{thm}

Theorem~\ref{mcmcthm} is proved in
Section~\ref{sec-mcmcproof} below.
Clearly, similar reasoning also applies with alternative assumptions
and to other versions
of adaptive MCMC including, for example, adaptive Metropolis-within-Gibbs
algorithms (with $P$ replaced by $P^d$ for random-scan);
cf. \cite{adapGibbs}.
\medskip

Theorem~\ref{mcmcthm} requires many conditions, but they are all
easy to ensure in practice, as illustrated
by the following type of adaptive MCMC algorithm:

\subsubsection*{The bounded adaption Metropolis (BAM) algorithm}
Let $\mathcal{X}= \mathbf{R}^d$, let $K \subseteq\mathcal{X}$ be
bounded, let $\pi$
be a continuous positive density on $\mathcal{X}$ and let $D>0$.
Let $\mathcal{Y}$ be a compact collection of $d$-dimensional
positive--definite matrices, and let $\Sigma_*\in\mathcal{Y}$ be fixed.
Define a process $\{X_n\}$ as follows: $X_0=x_0$ for some fixed $x_0\in K$.
Then for $n=0,1,2,\ldots,$ given $X_n$, we generate a proposal
$Y_{n+1}$ by:
(a) if $X_n \notin K$, then $Y_{n+1} \sim N(X_n, \Sigma_*)$;
(b) if $X_n \in K$ with $\operatorname{dist}(X_n,K^c) > 1$,
then $Y_{n+1} \sim N(X_n, \Sigma_{n+1})$, where
the matrix $\Sigma_{n+1}\in\mathcal{Y}$ is selected in
some fashion, perhaps depending on $X_n$ and
on the chain's entire history;
(c) if $X_n \in K$ but $\operatorname{dist}(X_n,K^c) = u$ with $0 \le
u \le1$,
then $Y_{n+1} \sim(1-u) N(X_n, \Sigma_*) + u N(X_n, \Sigma_{n+1})$.
Once $Y_{n+1}$ is chosen, then if $|Y_{n+1}-X_n|>D$, the
proposal is rejected so $X_{n+1}=X_n$.
Otherwise, if $|Y_{n+1}-X_n| \le D$, then with probability
$\min[1, {\pi(Y_{n+1}) \over\pi(X_n)}]$\vspace*{1pt} the proposal is accepted
so $X_{n+1}=Y_{n+1}$, or with the remaining probability the
proposal is rejected so $X_{n+1}=X_n$.

%re7 #&#
\begin{rem}
In the above BAM algorithm, $q_\gamma(Y_{n+1},X_n) = q_\gamma(X_n,Y_{n+1})$,
so those factors cancel in the formula for the acceptance probability.
\end{rem}

%re8 #&#
\begin{rem}
One good choice for the proposal covariance matrix $\Sigma_{n+1}$ in
part~(b) of the BAM algorithm is $(2.38)^2 V_n/d$, where
$V_n$ is the empirical covariance matrix of $X_0,\ldots,X_n$
from the process's previous history (except restricted to
some compact set $\mathcal{Y}$), since that choice approximates the optimal
proposal covariance; see the discussion in Section~2 of
\cite{adaptex.tex}.
\end{rem}

%pr22 #&#
\begin{prop}
\label{BAMprop}
The above BAM algorithm satisfies containment~(\ref{containmenteqn}).
Hence, if the selection of the $\Sigma_n$
satisfies diminishing adaptation~(\ref{dimadeqn}),
then convergence (\ref{MCMCconveqn}) holds.
\end{prop}

\begin{pf}
The BAM algorithm satisfies all of the conditions of Theorem~\ref{mcmcthm}.
Indeed, bounded jumps (\ref{boundedPgamma})
and no adaption outside of $K$ (\ref{fixedPeqn}) are both immediate.
Here the fixed kernel $Q$ is
bounded above (\ref{Pboundedabove}) by the constant $M = (2\pi)^{-d/2}
|\Sigma_*|^{-1/2}$, and
the $\epsilon$--$\delta$ condition (\ref{deltepeqn}) holds by the
formula for $Q$ together with the continuity of the density $\pi$
(which guarantees that it is bounded above and below on any compact
rectangle $J$ containing the compact set $K_{2D}$). Furthermore the
continuous densities condition (A6) holds by construction.
Hence, the result follows from Theorem~\ref{mcmcthm}.
\end{pf}

% this algorithm satisfies Containment~\eqref{containmenteqn}.
% It follows that if the selection of the matrices $\Sigma$ are done in
% a way which satisfies Diminishing Adaptation~\eqref{dimadeqn},
% then convergence \eqref{MCMCconveqn} holds.

% subject to the continuity condition that as a function of $X_n$
% the matrix $\Sigma_{n+1}(X_n)$ converges to $\Sigma_0$ as $X_n$
% converges to the boundary of $K$.

%s11.3 #&#
\subsection{Proof of Theorem~\texorpdfstring{\protect\ref{mcmcthm}}{21}}
\label{sec-mcmcproof}

We begin with a result linking the boundedness property
(\ref{boundedinprob}) for $\{X_n\}$ with the containment condition
(\ref{containmenteqn}) for $\{M_\epsilon(X_n,\break \Gamma_n)\}$, as follows:

%pr23 #&#
\begin{prop} \label{contprop}%
Consider an adaptive MCMC algorithm as in Section~\ref{sec-mcmcsetup}.
Suppose (\ref{boundedinprob}) holds, and for each $n\in\mathbf{N}$
the mapping
$(x,\gamma) \mapsto\Delta(x,\gamma,n)
:= \|P_\gamma^n(x,\cdot) - \pi(\cdot)\|$
is continuous. Then the containment condition~(\ref{containmenteqn}) holds.
\end{prop}

\begin{pf}
Since each $P_\gamma$ is Harris ergodic,
$\lim_{n\to\infty} \Delta(x,\gamma,n) = 0$
for each fixed $x\in\mathcal{X}$ and $\gamma\in\mathcal{Y}$.
Also, since $\pi$ is a stationary distribution for $P_\gamma$,
the mapping $n \mapsto\Delta(x,\gamma,n)$ is nonincreasing; see, for
example, Proposition~3(c) of \cite{probsurv.tex}.
If the mapping $(x,\gamma) \mapsto\Delta(x,\gamma,n)$ is continuous,
then it follows by Dini's theorem (e.g., \cite{rudin}, page~150) that
for any compact subset $C \subseteq\mathcal{X}$, since $\mathcal{Y}$
is compact,
\[
\lim_{n\to\infty} \sup_{x\in C} \sup
_{\gamma\in\mathcal{Y}} \Delta(x,\gamma,n) = 0.
\]
Hence, given $C$ and $\epsilon>0$, there is $n\in\mathbf{N}$ with
$\sup_{x\in C} \sup_{\gamma\in\mathcal{Y}} \Delta(x,\gamma,n) <
\epsilon$. It follows that
$\sup_{x\in C} \sup_{\gamma\in\mathcal{Y}} M_\epsilon(x,\gamma)
< \infty$
for any fixed $\epsilon>0$.

Now, if $\{X_n\}$ is bounded in probability as in (\ref{boundedinprob}),
then for any $\delta>0$,
we can find a large enough compact subset $C$ such that $P(X_n
\notin C) \le\delta$ for all $n$. Then given $\epsilon>0$, and
if $L := \sup_{x\in C} \sup_{\gamma\in\mathcal{Y}} M_\epsilon
(x,\gamma)$,
then $L<\infty$, and $P(M_\epsilon(X_n,\Gamma_n) > L) \le\delta$
for all $n$ as well.
Since $\delta$ was arbitrary, it follows that
$\{M_{\epsilon}(X_n,\Gamma_n)\}_{n=0}^\infty$ is bounded in probability.
\end{pf}

We then need a lemma guaranteeing continuity of
$\Delta(x,\gamma,n)$:

%le24 #&#
\begin{lem} \label{MHcont}
Under the continuous density assumptions \textup{(A6)},
for each $n\in\mathbf{N}$,
the mapping $(x,\gamma) \mapsto\Delta(x,\gamma,n)$ is continuous.
\end{lem}

\begin{pf}
Assuming (A6)(ii),
this fact is contained in the proof of Corollary~11
of \cite{adapt.tex}.
The corresponding result assuming (A6)(i)
is similar but easier.
\end{pf}

\begin{pf*}{Proof of Theorem~\protect\ref{mcmcthm}}
The bounded jumps condition (\ref{boundedPgamma}),
together with no adaption outside of $K$ (\ref{fixedPeqn}),
ensure that the algorithm $\{X_n\}$ fits within
the setup of Section~\ref{sec-setup}.
Since the densities of $P(x,dy)$ are bounded above by
(\ref{Pboundedabove}), it follows that (A1) holds
with $\mu_* = \operatorname{Uniform}(K_{2D} \setminus K_D)$.
Also, using the $\epsilon$--$\delta$ condition (\ref{deltepeqn}),
it follows from Proposition~\ref{unifsmallprop} that
(A3) holds for $\nu_* = \mu_*$.
Hence, by Theorem~\ref{thm-pismall}(a), $\{X_n\}$ is bounded
in probability; that is, (\ref{boundedinprob}) holds.
In addition, using the continuity assumption (A6),
it follows from Lemma~\ref{MHcont}
that $\Delta(x,\gamma,n)$ is a continuous function.
Containment~(\ref{containmenteqn})
thus follows from Proposition~\ref{contprop}.
The final assertion about convergence~(\ref{MCMCconveqn})
then follows from \cite{adapt.tex}.
\end{pf*}

%s12 #&#
\section{A detailed statistical MCMC example: RCA}
\label{sec-statexample}

Relying on the theoretical advances in this paper, we shall now
demonstrate the effectiveness of a general adaptive strategy which
we call \textit{regime change algorithm (RCA)} that can be implemented
in a wide number of practical instances. Specifically, during the
initialization period the chain is run using a transition kernel
that can provide some information about the target. We do not assume
that this initial kernel is optimal in any way, just that it would
be a reasonable initial choice for an MCMC algorithm. After the
initialization period, inside a chosen compact set, the initial kernel
is slowly replaced by an adaptive kernel that is shown to exhibit
better mixing. In a statistical example below, we shall see that
this regime change dramatically increases the algorithm efficiency,
since the adaptive kernel is increasingly more suitable for sampling
the target inside the compact. Our regime change idea is in the same
general vein as the two-stage adaptation proposed by Giordani and
Kohn \cite{kohn}. However, their theoretical justification follows
a rather different approach from ours.

%s12.1 #&#
\subsection{Model and data}
\label{sec-lupusdata}

%
%t1 #&#
\begin{table}[b]
\tablewidth=250pt
\caption{The number of latent membranous lupus nephritis
cases (numerator), and the total number of cases (denominator),
for each combination of the values of the two covariates,
for the 55 lupus patients
in the data set described in Section~\protect\ref{sec-lupusdata}}
\label{tb:dta}
\begin{tabular*}{250pt}{@{\extracolsep{\fill}}lccccc@{}}
\hline
& \multicolumn{5}{c@{}}{$\bolds{IgA}$}\\[-4pt]
& \multicolumn{5}{c@{}}{\hrulefill}\\
$\bolds{\Delta IgG}$& \textbf{0} & $\bolds{0.5}$ &$\bolds
{1}$ &$\bolds{1.5}$ &$\bolds{2}$\\
\hline
$-3.0$& $0/1$ &--&--&--&--\\
$-2.5$& $0/3$&--&--&--&--\\
$-2.0$& $0/7$&--&--&--&$0/1$\\
$-1.5$& $0/6$ & $0/1$ &--&--&--\\
$-1.0$& $0/6$ & $0/1$ & $0/1$ &--& $0/1$\\
$-0.5$& $0/4$&--&--& $1/1$ &--\\
\phantom{$-$}0&$0/3$&--& $0/1$ & $1/1$ &--\\
\phantom{$-$}0.5& $3/4$ &--& $1/1$ & $1/1$& $1/1$\\
\phantom{$-$}1.0& $1/1$ &--& $1/1$ & $1/1$ & $4/4$\\
\phantom{$-$}1.5& $1/1$ &--&--& $2/2$ &--\\
\hline
\end{tabular*}
\end{table}

We shall consider a Bayesian probit regression model applied
to a well-known collection of lupus patient data originally supplied
by Haas~\cite{haas} and later simplified in \cite{vanDyk}.
The data, shown in Table~\ref{tb:dta}, contain disease status for 55 patients
of which 18 have been diagnosed with latent membranous lupus, together
with two clinical covariates, $IgA$ and $\Delta IgG$ (which is equal
to $IgG3-IgG4$ in the lupus context), which are
computed from their levels of immunoglobulin of type $A$ and of type
$G$, respectively.
We consider a probit regression (PR) model; that is,
for each patient $1\le i\le55$, and we model the disease indicator
variables as independent
\[
Y_i \sim\operatorname{Bernoulli}\bigl(\Phi\bigl(x_i^T
\beta\bigr)\bigr), % \eqlabel{pr}
\]
where $\Phi(\cdot)$ is the CDF of
$N(0,1)$, $x_i=(1,\Delta IgG_{i},IgA_{i})$ is the vector of covariates,
and $\beta$ is a
$3\times1$ vector of parameters which is assigned
a flat prior $p(\beta)\propto1$.
The posterior is thus
\begin{eqnarray*}
&& \pi_{\mathrm{PR}}(\vec\beta|\vec Y, \vec IgA, \vec\Delta IgG)
\\
&& \qquad\propto \prod_{i=1}^{55} \bigl[
\Phi(\beta_{0} + \Delta IgG_{i}\beta_{1} +
IgA_{i}\beta_{2})^{Y_{i}}\\
&&\hspace*{6pt}\qquad\qquad{}\times\bigl(1-\Phi(
\beta_{0} \Delta IgG_{i}\beta_{1}+IgA_{i}
\beta_{2})\bigr)^{(1-Y_{i})} \bigr]. % \eqlabel{post_pr}
\end{eqnarray*}
We wish to design effective algorithms to sample from
this posterior distribution~$\pi_{\mathrm{PR}}$.

%s12.2 #&#
\subsection{The best previous algorithm: PX-DA}

The current state-of-the-art most efficient algorithm to sample from
the above posterior distribution $\pi_{\mathrm{PR}}$ is the \textit{parameter
expanded data augmentation (PX-DA)} algorithm developed by van Dyk
and Meng \cite{vanDyk}. The PX-DA transition kernel for updating
$\beta^{(t)}$ is defined by the following steps:
\begin{itemize}
\item Draw
\[
\phi_{i}^{(t+1)}\sim\cases{ N_{+}
\bigl(x_{i}^{T}\beta^{(t)},1 \bigr),& \quad$
\mbox{if }Y_{i}=1,$ \vspace*{3pt}
\cr
N_{-}
\bigl(x_{i}^{T}\beta^{(t)},1 \bigr),&\quad$
\mbox{if } Y_{i}=0,$}
\]
where $N_+(\mu,\sigma^{2})$ and $N_-(\mu,\sigma^{2})$ are
normal distributions with mean $\mu$ and variance $\sigma^{2}$
that are truncated to $(0,\infty)$ and $(-\infty,0)$, respectively.
Set $\phi^{(t+1)}=(\phi_{1}^{(t+1)},\ldots,\phi_{n}^{(t+1)})$.

\item Let $\tilde\beta^{t+1}= (X^{T}X)^{-1}X^{T}\phi^{(t+1)}$, and define\vspace*{2pt}
$R^{(t+1)}=\sum_{i=1}^{n}(\phi_{i}^{(t+1)}- \break x_{i}^{T}\tilde\beta
^{(t+1)})^{2}$.

\item Sample $Z \sim N(0,1)$, $W\sim\chi_{n}^{2}$ and set
$\beta^{(t+1)}=\sqrt{{W \over R^{(t+1)}}}\tilde\beta^{(t+1)}
+ \break \operatorname{Chol}[(X^{T} X)^{-1}]Z$.
\end{itemize}

%s12.3 #&#
\subsection{A new algorithm: RCA}

The regime change algorithm (RCA) is initialized by running the PX-DA
chain for $M$ iterations. Based on the samples obtained, we determine a
compact subset $K$ and a distance bound $D$
which remain fixed for the rest of the simulation.
The algorithm then proceeds by constructing a Gaussian approximation
of the target inside $K$ that continuously evolves as the samples are
collected, thus allowing for better and better proposal values.

To proceed, for $n \ge M$ we define
\[
\mu_{n} := \frac{ \langle X_0 \rangle + \langle X_1 \rangle +
\cdots+ \langle X_{n-1} \rangle }{n},
\]
and
\[
\Sigma_{n} := \operatorname{Cov}\bigl(\langle X_0
\rangle ,\langle X_1 \rangle ,\ldots,\langle X_{n-1} \rangle
\bigr) + \epsilon I_d,
\]
where $\operatorname{Cov}$ is the empirical covariance function, and
$\langle r \rangle $ is the shrunken version of $r\in\mathbf{R}^d$
with each coordinate
shrunk into the interval $[-L,L]$, that is,
$\langle r \rangle _i = \max[-L, \min(L,r_i)]$.
We then define $K$ to be the ball centered at $\mu_{M}$,
of radius $\max_{1 \le i \le d} (\Sigma_M)_{ii}^{1/2}$
(i.e., the largest sample standard
deviation on the diagonal of $\Sigma_M$), and we let $D$ be any
suitably large distance bound (e.g., $D=20$).

We then consider the independence Metropolis (IM) transition kernel
$P_{\mu_n,\Sigma_n}$, with proposal distribution given (independently
of the current state of the process) by the Gaussian distribution
$N(\mu_{n},\Sigma_{n})$, except truncated (in a continuous manner;
see Remark~\ref{contremark} below) to remain in the compact $K$ and
to never move more than a distance $D$. We also let $P_{PX}(x,y)$ be
the PX-DA algorithm described above, also truncated in a continuous
manner to remain in the compact $K$ and to never move more than a
distance $D$.

In terms of these definitions, the update for the RCA follows these
steps:
\begin{longlist}[(3)]
\item[(1)] If $X_n \in K^c$, then $X_{n+1} \sim P_{PX}(X_n,\cdot)$.

\item[(2)] If $X_n \in K$ and $d(X_n,K^c) > 1$, then
\[
X_{n+1} \sim \lambda_{n+1} P_{\mu_n,\Sigma_n}(X_n,
\cdot) + (1-\lambda_{n+1}) P_{PX}(X_n, \cdot),
\]
with $\lambda_{n} = \min[ \max( \theta_{n},0.2), 0.8]$,
where $\theta_{n}$ is the empirical acceptance rate of all of
the IM proposals made so far between time $M+1$ and time $n-1$
(or we simply set $\lambda_n=1/2$ if there have been no such proposals).

\item[(3)] If $X_n \in K$ and $d(X_n,K^c)=u$ with $0 \le u \le1$, then
\[
X_{n+1} \sim u \bigl[ \lambda_{n+1} P_{\mu_n,\Sigma_n}(X_n,
\cdot) + (1-\lambda_{n+1}) P_{PX}(X_n, \cdot)
\bigr] + (1-u) P_{PX}(X_n, \cdot),
\]
with $\lambda_n$ as above.
\end{longlist}

That is, letting $\gamma_n=(\mu_n,\Sigma_n,\lambda_n)$ be the complete
adaptive parameter, we can say that when
$d(X_{n},K^{c}) > 1$, the chain moves according to the adaptive kernel
\[
P_{K,\gamma_{n}}(X_{n},\cdot) = \lambda_{n+1}
P_{\mu_n,\Sigma_n}(X_n,\cdot) + (1-\lambda_{n+1})
P_{PX}(X_n, \cdot),
\]
and when $X_{n}\in K^{c}$ the chain follows the transition
$P_{PX}(X_{n},\cdot)$,
with a linear interpolation near the boundary of $K$ to satisfy
the continuous densities condition~(A6).

%re9 #&#
\begin{rem}
\label{contremark}
In our description of RCA above, we required certain Gaussian distributions
to be restricted to certain subsets.
If this is done naively, then it will result
in a discontinuous density, which may violate (A6).
However, this issue can be easily avoided if we make
the density continuous by smoothing the edge via
a linear interpolation. For example, to restrict a
univariate normal density with mean $\mu$ and variance $\sigma^2$
to the range $(a,b)$ for $a<b$, one can choose small $\upsilon>0$
and define
\[
f_\upsilon(x | \mu,\sigma,a,b) = \frac{(2\pi\sigma^2)^{-1/2} \exp [-{(x-\mu)^2}/(2\sigma^2)]
}{\Phi(b-\upsilon)-\Phi(a+\upsilon)},
\]
and then use the density function proportional to
\[
g(x | \mu,\sigma,a,b,\upsilon)= \cases{ f_\upsilon(x | \mu,\sigma,a,b), &
\quad$\mbox{if } a+\upsilon\le x \le b-\upsilon$,\vspace*{3pt}
\cr
f_\upsilon(b-\upsilon | \mu,\sigma,a,b) (b-x)/\upsilon, & \quad$\mbox{if }
b-\upsilon< x < b,$\vspace*{3pt}
\cr
f_\upsilon(a+\upsilon | \mu,
\sigma,a,b) (x-a)/\upsilon, & \quad$\mbox{if } a< x <a+\upsilon$, \vspace*{3pt}
\cr
0, & \quad$\mbox{otherwise}$.}
\]
The general multivariate case can be handled by similarly truncating
each of the independent univariate Gaussian variables used to
construct the multivariate Gaussian. In this way, it can be assured
that even truncated Gaussians still have continuous densities.
\end{rem}

%s12.4 #&#
\subsection{Verification of the theoretical assumptions}

To justify the use of our new RCA algorithm, we wish to prove
asymptotic convergence as in (\ref{MCMCconveqn}). Proving such
convergence of adaptive MCMC algorithms is usually very difficult, but
we shall manage this by applying Theorem~\ref{mcmcthm}. To do this,
we need to verify the assumptions of Theorem~\ref{mcmcthm} including
those which are implicit in the setup of Section~\ref{sec-mcmcsetup}.
Fortunately, this is not too difficult.

For the RCA algorithm, the ``bounded jumps'' condition
(\ref{boundedPgamma}) and the ``fixed kernel outside of $K$''
condition (\ref{fixedPeqn}) are both satisfied by construction.

Furthermore, the ``fixed kernel bounded above by a multiple
of Lebesgue'' condition (\ref{Pboundedabove}), and the
``$\epsilon$--$\delta$ bounded below by a multiple of Lebesgue''
condition~(\ref{deltepeqn}), both concern the transition probabilities
outside of $K$, and hence they both follow since our fixed transition
probabilities are absolutely continuous with respect to Lebesgue
measure with densities that are uniformly bounded away from 0 and
$\infty$ on compact subsets.

In addition, the continuous densities condition (A6) is
satisfied since all transition kernels involved in the construction of
the chain are Metropolis--Hastings (MH) kernels with proposal densities
that are continuous functions of the adaption parameters and of
$x$ and $y$; cf. Remark~\ref{contremark}.

Finally, we note that RCA also satisfies the diminishing adaptation
condition~(\ref{dimadeqn}), since the difference between the values
of each of the adaptation parameters at iterations $n$ and $n+1$
is always $O(n^{-1})$.

Hence RCA satisfies all of the assumptions of Theorem~\ref{mcmcthm}
and Section~\ref{sec-mcmcsetup}, and also
satisfies Diminishing Adaptation~(\ref{dimadeqn}),
so we conclude:

%co25 #&#
\begin{cor}
\label{RCAcor}
The RCA algorithm described above converges asymptotically to $\pi$
as in~(\ref{MCMCconveqn}).
\end{cor}

%s12.5 #&#
\subsection{A simulation study}
\label{sec-simulation}

To test our new RCA algorithm in practice, we ran\footnote{The R
computer program we used is available at:
\surl{www.probability.ca/lupus}}
both it and the PX-DA
algorithm, each for 5000 iterations starting with $X_0$ equal to
the maximum likelihood estimate (MLE).

\begin{figure}

\includegraphics{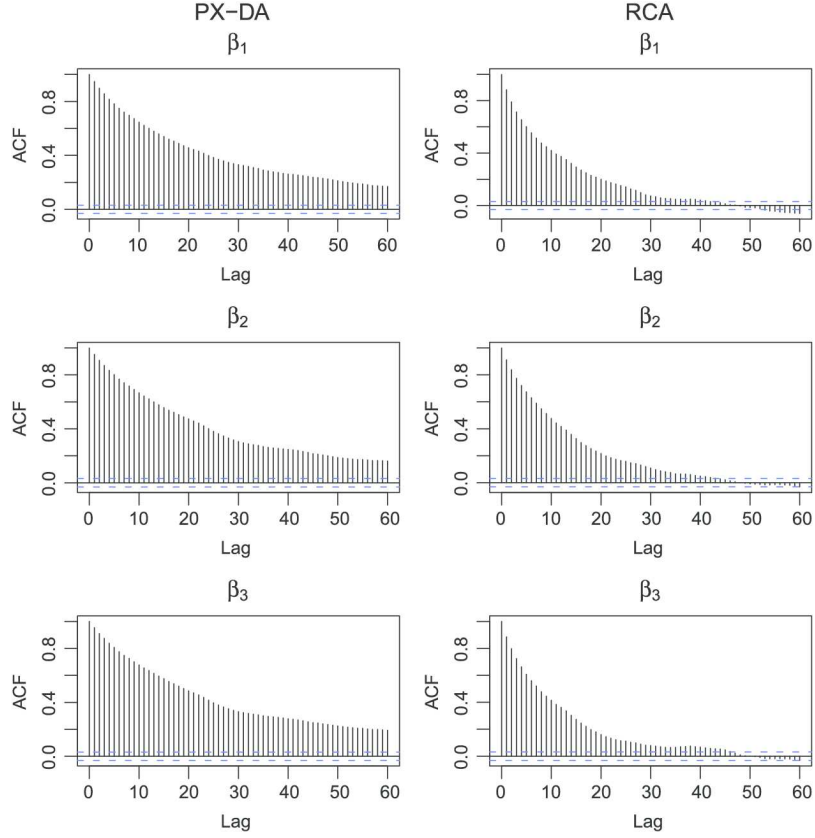}

\caption{Autocorrelation (ACF) plots for the
probit regression model simulation study % \hfil\break
of Section~\protect\ref{sec-simulation}, comparing the PX-DA (left column)
and RCA (right column) algorithms, % \hfil\break
for each of the three parameters $\beta_0$ (top),
$\beta_1$ (middle) and $\beta_2$ (bottom), showing significantly
smaller autocorrelations (and hence better performance)
for RCA than for PX-DA.}
\label{fig:lupus_acf}
\end{figure}

We found that the RCA algorithm did indeed perform significantly more
efficiently than PX-DA did. As one measure of this, we plotted the
autocorrelation function (ACF) plots of both algorithms for each of
the three parameters (Figure~\ref{fig:lupus_acf}). This plot indicates
that the autocorrelations for RCA are significantly smaller than
those for PX-DA, thus indicating faster mixing and thus a more
efficient algorithm. Indeed, the sums of the nonnegligible
positive-lag autocorrelations for the three parameters were respectively
41.20, 40.87 and 43.87 for PX-DA, but just
10.56, 11.88 and 10.00 for RCA, and again showing much greater
efficiency of RCA.

%
%f3 #&#

Another way to think about this is in terms of
effective sample size (ESS). This is a measure of how many true
independent samples our algorithm is equivalent to, in terms of
variance of the resulting estimator. The ESS is well known
(see, e.g., \cite{flegal}, page~2) to be inversely proportional to
$1 + 2S$ where $S$ is the autocorrelation sum as above. By this
measure, in our simulations the ESS for RCA is larger than for PX-DA,
for the three parameters respectively, by factors of
$3.77$, $3.34$ and $4.23$.
This indicates quite significant improvements in efficiency
of RCA over PX-DA for this example.

% $(1+2*41.20)/(1+2*10.56) = 3.77$,
% $(1+2*40.87)/(1+2*11.88) = 3.34$,
% and $(1+2*43.87)/(1+2*10.00) = 4.23$.

We conclude that having the possibility to sample from the IM
kernel reduces the autocorrelation within the samples produced
by the algorithm and thus significantly increases the effective
sample size. This indicates that the RCA algorithm (as justified in
Corollary~\ref{RCAcor}, by applying Theorem~\ref{mcmcthm}) is indeed
a superior algorithm for this problem.

%sA #&#
\begin{appendix}
%sB #&#
\section*{Appendix: Replacing the minorizing measure by $\pi$}
\label{sec-appendix}

Recall that assumption~(A3) requires that the set $K_{2D}
\setminus K_D$ be small for $P$, with some minorizing measure $\nu_*$.
It turns out that if
assumption~(A3) holds for any $\nu_*$, and if
$P$ is reversible, then assumption~(A3) also holds
for the specific choice $\nu_* = \pi|_{K_{2D} \setminus K_D}$, that is,
where $\nu_*(A) = \pi(A \cap(K_{2D} \setminus K_D)) /
\pi(K_{2D} \setminus K_D)$,
with the step size $n_0$ replaced by $2n_0$.
Under the additional assumption of uniform ergodicity,
this fact is Proposition~1 of \cite{sechyb.tex}.
For arbitrary reversible chains, this fact follows from
Lemma~5.9 of the Polish doctoral thesis \cite{KrysBrother}, which for
completeness we now reproduce:

%le26 #&#
\begin{lem}[(Lemma 5.9 of \cite{KrysBrother})]
\label{BM_lemma}
Let $P$ be a Markov chain transition kernel on $(\mathcal{X},\mathcal{F})$,
with invariant probability measure $\pi$.
Let $C\in\mathcal{F}$ such that $\pi(C)>0$.
Assume that $C$ is a small set for
$P$; that is, for some $n_0\in\mathbf{N}$ and $\beta>0$ and
probability measure $\nu$,
%
%eB.1 #&#
\begin{equation}
\label{general_small} P^{n_0}(x,A) \ge \beta \mathbf{1}_C(x)
\nu(A),\qquad A \in\mathcal{F}.
\end{equation}
Then
%
%eB.2 #&#
\begin{equation}
\label{pi_small} P^{n_0} \bigl(P^*\bigr)^{n_0}(x,A) \ge
\tfrac{1}{4} \beta^2 \mathbf{1}_C(x) \pi(A\cap C),
\qquad A \in\mathcal{F},
\end{equation}
where $P^*$ is the $L^2(\pi)$ adjoint of $P$.
In particular, if $P$ is reversible with respect to $\pi$, so
that $P^*=P$, then
\[
P^{2n_0}(x,A) \ge\tfrac{1}{4} \beta^2
\mathbf{1}_C(x) \pi(A\cap C),\qquad A \in\mathcal{F}.
\]
\end{lem}

Hence if $K_{2D} \setminus K_D$ is an $n_0$-small set with minorizing
measure $\nu$, and $P$ is reversible with respect to $\pi$,
then $K_{2D} \setminus K_D$
is a $(2 n_0)$-small set with minorizing measure $\pi|_{K_{2D}
\setminus K_D}$.

\begin{pf*}{Proof of Lemma~\ref{BM_lemma}}
By replacing\vspace*{1pt} $P$ by $P^{n_0}$ and $P^*$ by $(P^*)^{n_0}$,
it suffices to assume that $n_0=1$.
Now, the Radon--Nikodym derivative $\frac{d\nu}{d\pi}$
of $\nu$ with respect to $\pi$ satisfies that
%Integrating both sides of~\eqref{general_small}
%with respect to $\pi$ yields that
%$$
%\frac{d\nu}{d\pi} \ \le\ \beta\pi(C)
% {\rm and}
$\int_{\mathcal{X}}\frac{d\nu}{d\pi}(x) \pi(dx) = \nu(\mathcal
{X}) = 1$.
Hence, for every $\varepsilon\in[0,1]$, the set
%
%eB.3 #&#
\begin{equation}
\label{def_D} D(\varepsilon) := \biggl\{ x\in\mathcal{X} \dvtx
\frac{d\nu}{d\pi}(x) \ge \varepsilon \biggr\}
\end{equation}
has $\pi(D(\varepsilon)) > 0$. We then compute that
\[
\nu\bigl(D(\varepsilon)^c\bigr) = \int_{D(\varepsilon)^c}
\frac{d\nu}{d\pi}(x) \pi(dx) \leq \varepsilon \int_{\mathcal{X}}
\pi(dx) = \varepsilon %
\]
and hence
%
%eB.4 #&#
\begin{equation}
\label{bound_D} \nu\bigl(D(\varepsilon)\bigr) \ge 1-\varepsilon.
\end{equation}
Recall also that the adjoint $P^*$ satisfies
%
%eB.5 #&#
\begin{equation}
\label{p_star} \pi(dx) P(x, dy) = \pi(dy) P^*(y,dx).
\end{equation}
Now let $x \in C$, and $A\in\mathcal{F}$ with $A\cap C \neq\varnothing$.
Using first (\ref{general_small}) and then (\ref{def_D}),
\begin{eqnarray*}
PP^*(x, A) &=& \int_{z \in\mathcal{X}} P^*(z, A) P(x, dz) \ge \beta\int
_{z \in\mathcal{X}} P^*(z, A\cap C) \nu(dz)
\\
&\ge& \beta\int_{z \in D(\varepsilon)} \int_{y \in A \cap C}
P^*(z, dy) \varepsilon \pi(dz).
\end{eqnarray*}
To continue, use~(\ref{p_star}), then (\ref{general_small}) again
and finally (\ref{bound_D}) to obtain
\begin{eqnarray*}
PP^*(x, A) &\ge& \beta \varepsilon \int_{z \in D(\varepsilon)} \int
_{y \in A \cap C} \pi(dy) P(y, dz)
\\
&\ge& \beta^2 \varepsilon \nu\bigl(D(\varepsilon)\bigr) \pi(A\cap C)
\ge \beta^2 \varepsilon(1-\varepsilon) \pi(A\cap C).
\end{eqnarray*}
Setting $\varepsilon=1/2$ yields (\ref{pi_small}).
\end{pf*}
\end{appendix}

\section*{Acknowledgements}
We thank Blazej Miasojedow and Daniel Rudolf for helpful comments,
and thank the two anonymous referees for very careful readings of
the paper which led to significant improvements.

% imsref loaded by daiva.urboniene, 2015-01-06 16:52:48
% imsref loaded by daiva.urboniene, 2015-01-06 17:20:02
%

%\begin{appendix}
%\section{}
%\end{appendix}

%\section*{Acknowledgments}

%\begin{supplement}[id=suppA]
%\slink[doi]{10.1214/00-AAPXXXXSUPP} %[doi,text={...}] - jei reikia
%suskaldyti doi
%\sfilename{aapXXXX\_supp.pdf}
%\end{supplement}

%\begin{thebibliography}{99}
%\bibitem[\protect\citeauthoryear{}{}]{r1}
%\bibitem{r1}
%\end{thebibliography}

\printaddresses
\end{document}